\documentclass[nonblindrev]{informs3_noinforms}
\usepackage{pdfsync}
\usepackage{enumerate}

\DoubleSpacedXI

\usepackage{endnotes}
\let\footnote=\endnote

\newcommand{\E}{\bbe}

\newcommand{\p}{\mathbb{P}}

\def\bbr{{\Bbb{R}}} 
\def\bbe{{\Bbb{E}}} 

\renewcommand{\theequation}{\thesection.\arabic{equation}}

\usepackage{natbib}
 \bibpunct[, ]{(}{)}{,}{a}{}{,}%
 \def\bibfont{\small}%
\TheoremsNumberedThrough     
\ECRepeatTheorems

\EquationsNumberedThrough    

\begin{document}



\RUNTITLE{Asymptotic optimality of TBS policies in dual-sourcing inventory systems}
\TITLE{Asymptotic optimality of Tailored Base-Surge policies in dual-sourcing inventory systems}

\ARTICLEAUTHORS{
\AUTHOR{\bf Linwei Xin}
\AFF{Department of Industrial and Enterprise Systems Engineering,
University of Illinois at Urbana-Champaign,
Urbana, IL 61801}
\AUTHOR{\bf David A. Goldberg}
\AFF{School of Industrial \& Systems Engineering,
Georgia Institute of Technology,
Atlanta, GA 30332}

\RUNAUTHOR{{Xin and Goldberg}}

} 

\ABSTRACT{\indent Dual-sourcing inventory systems, in which one supplier is faster (i.e. express) and more costly, while the other is slower (i.e. regular) and cheaper, arise naturally in many real-world supply chains.  These systems are
notoriously difficult to optimize due to the complex structure of the optimal solution and the curse of dimensionality, having resisted solution for over 40 years.  Recently, so-called Tailored Base-Surge (TBS) policies have been proposed as a heuristic for the dual-sourcing problem.  Under such a policy, a constant order is placed at the regular source in each period, while the order placed at the express source follows a simple order-up-to rule.  Numerical experiments by several authors have suggested that such policies perform well as the lead time difference between the two sources grows large, which is exactly the setting in which the curse of dimensionality leads to the problem becoming intractable.  However, providing a theoretical foundation for this phenomenon has remained a major open problem.
\\\indent In this paper, we provide such a theoretical foundation by proving that a simple TBS policy is indeed asymptotically optimal as the lead time of the regular source grows large, with the lead time of the express source held fixed.  Our main proof technique combines novel convexity and lower-bounding arguments, an explicit implementation of the vanishing discount factor approach to analyzing infinite-horizon Markov decision processes, and ideas from the theory of random walks and queues, significantly extending the methodology and applicability of a novel framework for analyzing inventory models with large lead times recently introduced by Goldberg and co-authors in the context of lost-sales models with positive lead times.
}



\maketitle

%


\section{Introduction}
A common practice in the management of global supply chains is \emph{dual-sourcing} (cf. \citet{Rao00}).  Under a dual-sourcing strategy, the companies usually purchase their materials from a regular supplier at a lower cost, but they are also able to obtain materials from an expedited supplier at a higher cost under emergency circumstances. For example, in the summer of 2003, Amazon used FedEx to deliver the new Harry Potter more promptly and maintained regular shipping via UPS (cf. \citet{Kelleher03}, \citet{Veeraraghavan08}). \citet{Allon10} describes an example of a $\$10$ billion high-tech U.S. company that has two suppliers, one in Mexico and one in China. The one in Mexico has shorter lead time but higher per-unit ordering cost; the one in China has longer lead time (5 to 10 times longer) but lower per-unit ordering cost. The company takes advantage of the dual-sourcing strategy to meet the demand more responsively (from Mexico) as well as less expensively (from China).

Although dual-sourcing is attractive, and very relevant to practice, optimizing a dual-sourcing inventory system is notoriously challenging.  Such inventory systems have been studied now for over forty years, but the structure of the optimal policy remains poorly understood, with the exception of when the system is consecutive, i.e., the lead time difference between the two sources is exactly one.  More specifically, the earliest studies of periodic review dual-sourcing inventory models include \citet{Barankin61}, \citet{Daniel63}, and \citet{Neuts64}, which showed that base-stock (also known as order-up-to) policies are optimal when the lead times of the two sources are zero and one respectively. \citet{Fukuda64} extended the result to general lead time settings as long as the lead time difference remains one. \citet{Whittmore77} showed that the optimal policy is no longer a simple base-stock policy when the lead time difference is beyond one and the structure of the optimal policy can be quite complex. Furthermore, it is well known that a dual-sourcing inventory system can be regarded as a generalization of a lost-sales inventory system (cf. \citet{Sheopuri10}).  Indeed, the intractability of both the dual-sourcing and lost-sales inventory models has a common source - as the lead time grows, the state-space of the natural dynamic programming (DP) formulation grows exponentially, rendering such techniques impractical.  This issue is typically referred to as the ``curse of dimensionality" (cf. \citet{Karlin58}, \citet{Morton69}, \citet{Zipkin08b}), and we refer the reader to \citet{Goldberg12} and \citet{Xin14a} for a relevant discussion in the context of lost-sales inventory models.

There is a vast literature investigating periodic review dual-sourcing inventory models as well as their variants, and we refer the interested reader to the survey of \citet{Minner03}, as well as e.g. the more recent works of \citet{Feng06}, \citet{Fox06}, \citet{Chen13}, \citet{Huggins10}, \citet{Angelus15}, \citet{Boute14}, \citet{Gong14}, \citet{Song09}, and the references therein.

As an exact solution seems out of reach, the operations research and management communities have instead investigated certain structural properties of the optimal policy (cf. \citet{Hua14}), and exerted considerable effort towards constructing various heuristic policies. \citet{Veeraraghavan08} proposed the family of dual index (DI) policies, which have two base-stock levels, one for the regular source and one for the express source, and ``orders up" to bring appropriate notions of inventory position up to these levels.  \citet{Scheller-Wolf08} analyzed the closely related class of single index (SI) policies, for which the relevant notions of inventory position are different.  Both families of policies seem to perform well in numerical studies.
\citet{Sheopuri10} considered two generalized classes of policies: one with an order-up-to structure for the express source, and one with an order-up-to structure for the regular source. Their numerical experiment showed that such policies can outperform DI policies.  In the presence of production capacity costs, \citet{Boute14} studied dual-sourcing smoothing policies, under which the order quantities from both sources in each period are convex combinations of observed past demands. They analyzed such polices under normally distributed demand, and their numerical results showed that these policies performed better for higher capacity costs and longer lead time differences (between the two sources).

A simple and natural policy that is implemented in practice, which will be the subject of our own investigations, is the so-called \emph{Tailored Base-Surge} (TBS) policy. It was first proposed and analyzed in \citet{Allon10}, where we note that closely related standing order policies had been studied previously (cf. \citet{Rosenshine76,Janssen99}).  Under such a TBS policy, a constant order is placed at the regular source in each period to meet a \emph{base} level of demand, while the orders placed at the express source follow an order-up-to rule to manage demand \emph{surges}.  We refer to Mini-Case 6 in \citet{Van Mieghem08} for more about the motivation and background of TBS policies. Note that dual-sourcing inventory systems in which a constant-order policy is implemented for the regular source are essentially equivalent to single-sourcing inventory systems with constant returns, which have been investigated in the literature (cf. \citet{Fleischmann03}, \citet{DeCroix05}).

\citet{Allon10} analyzed TBS policies in a continuous review model, and their focus was to find the best TBS policy. Numerical results in \citet{Klosterhalfen11} and \citet{Rossi12} showed that TBS policies are comparable to DI policies, and outperform DI policies for some problem instances. \citet{Allon10} conjectured that this policy performs more effectively as the lead time difference between the two sources grows. \citet{Janakiraman14a} (henceforth denoted JSS) analyzed a periodic review model and studied the performance of the TBS policy. They provided an explicit bound on the performance of TBS policies compared to the optimal one when the demand had a specific structure, and provided numerical experiments suggesting that the performance of the TBS policy improves as the lead time difference grows large.

However, to date there is no theoretical justification for the good behavior of TBS policies as the lead time difference grows large, and giving a solid theoretical foundation to this observed phenomena remains a major open question.  We note that until recently, a similar state of affairs existed regarding the good performance of constant-order policies as the lead time grows large in single-source lost-sales inventory models.  However, using tools from applied probability, queueing theory, and convexity, this phenomena was recently explained in \citet{Goldberg12} and \citet{Xin14a}, in which it was proven that a simple constant-order policy is asymptotically optimal in this setting as the lead time of the single source grows large.  The intuition here is that as the lead time grows large, so much randomness is introduced into the system between when an order is placed and when that order is received, that it is essentially impossible for any algorithm to meaningfully use the state information to make significantly better decisions.  Thus a policy which ignores the state information (i.e. constant-order policy) performs nearly as well as an optimal policy.  We note that the results of \citet{Xin14a} further demonstrate that the optimality gap of the constant-order policy actually shrinks exponentially fast to zero as the lead time grows large, and provide explicit and effective bounds even for moderate-to-small lead times.

\subsection{Our contributions}
In this paper, we resolve this open question by proving that, when the lead time of the express source is held fixed, a simple TBS policy is asymptotically optimal as the lead time of the regular source grows large.  Our results provide a solid theoretical foundation for the conjectures and numerical experiments of \citet{Allon10} and JSS.   Interestingly, the simple TBS policy performs nearly optimally exactly when standard DP-based methodologies become intractable due to the aforementioned ``curse of dimensionality".  Furthermore, as the ``best" TBS policy can be computed by solving a convex program that does not depend on the lead time of the regular source (cf. JSS), our results lead directly to very efficient algorithms (with complexity independent of the lead time of the regular source) with asymptotically optimal performance guarantees.  We also explicitly bound the optimality gap of the TBS policy for any fixed lead time (of the regular source), and prove that this decays inverse-polynomially in the lead time of the regular source.  Perhaps most importantly, since many companies are already implementing such TBS policies (cf. \citet{Allon10}), our results provide strong theoretical support for the widespread use of TBS policies in practice.  Our main proof technique combines novel convexity and lower-bounding arguments, 
an explicit implementation of the vanishing discount factor approach to analyzing infinite-horizon Markov decision processes (MDP), and ideas from the theory of random walks and queues.  Our methodology significantly extends the framework for analyzing inventory models with large lead times recently introduced in \citet{Goldberg12} and \citet{Xin14a} in the context of lost-sales models with positive lead times.  Indeed, in the present work we relate the performance of an optimal policy to a certain \emph{dynamic} optimization problem by applying the conditional Jensen's inequality, while in \cite{Xin14a} the relevant optimal policy could be bounded by a \emph{static} optimization problem after applying the (non-conditional) Jensen's inequality.  The inherently dynamic nature of the resulting bounds introduce several additional difficulties not encountered previously, and which we address in the present work.

\subsection{Outline of paper}
The rest of the paper is organized as follows. We formally define the dual-sourcing problem in Section \ref{ds-sec-intro}, and describe the TBS policy in Section \ref{ds-sec-TBS}.  We state our main result in Section \ref{ds-sec-main}, and prove our main result in  Section\ \ref{ds-sec-proof}.  We summarize our main contributions and propose directions for future research in Section \ref{ds-sec-conclusion}.  We also include a technical appendix in Section\ \ref{techsec}.

\section{Model description, problem statement and assumptions}\label{ds-sec-intro}
In this section, we formally define our dual-sourcing inventory problem, closely following the definitions given in \citet{Sheopuri10}.
Let $\{D_t\}_{t \in (-\infty,\infty)}, \{D'_t\}_{t \in (-\infty,\infty)}$ be mutually independent sequences of nonnegative independent and identically distributed (i.i.d.) demand realizations, distributed as the non-negative random variable (r.v.) $D$, which we assume to have finite mean, and (to rule out certain trivial degenerate cases) to have strictly positive (possibly infinite) variance.  Here we have introduced two doubly indexed sequences to prevent any possible confusion regarding dependencies of various demand realizations.  Let $\hat{G}$ be an independent geometrically distributed r.v., where $\p(\hat{G} = k) = 2^{-k}, k \geq 1$.  As a notational convenience, let us define all empty sums to equal zero, empty products to equal one, $\frac{1}{\infty} = 0$, $\mathbf{0} (\mathbf{1})$ denote the all zeros (ones) vector, and $\mathbb{I}(A)$ denote the indicator of the event $A$.  Let $L \geq 1$ be the deterministic lead time of the regular source (R), and $L_0 \geq 0$ the deterministic lead time of the express source (E), where $L > L_0 + 1$. Let $c_R,c_E$ be the unit purchase costs of the regular and express sources, and $h,b$ be the unit holding and backorder costs respectively, with $c \triangleq c_E - c_R > 0$.  In addition, let $I_t$ denote the on-hand inventory at the start of period $t$ (before any orders or demands are received), and $q_t^R (q_t^E)$ denote the order placed from R(E) at the beginning of period $t$.  Note that due to the leadtimes, the order received from R(E) in period $t$ is $q_{t - L}^R (q_{t - L_0}^E)$.  As we will be primarily interested in the corresponding long-run-average problem, and for simplicity (in later proofs), we suppose that the initial conditions are such that (s.t.) the initial inventory is $-\sum_{i=1}^{\hat{G}} D'_{-i}$, and no initial orders have been placed from either R or E.  Indeed, the associated system state will prove convenient to use as a ``regeneration point" when analyzing certain Markov chains which arise in our proofs, where we note that the geometric distribution allows us to preclude certain kinds of pathological periodic / lattice behavior which might otherwise interfere with proving the existence of relevant stationary measures.  We note that although assuming such a convenient randomized initial condition simplifies several technical proofs along these lines, such an assumption is not strictly necessary for our analysis, since the associated long-run average problem is insensitive to the particular choice of initial conditions.

As a notational convenience, we define $q_k^R = q_k^E = 0, k \leq 0$.  For $t = 1,\ldots,T$, the events in period $t$ are ordered as follows.
\begin{itemize}
\item Ordering decisions from R and E are made (i.e. $q_t^E,q_t^R$ are chosen);
\item New inventory $q_{t-L}^R + q_{t-L_0}^E$ is delivered and added to the on-hand inventory;
\item The demand $D_t$ is realized, costs for period $t$ are incurred, and the inventory is updated.
\end{itemize}
Note that the on-hand inventory is updated according to $I_{t+1}= I_{t} + q_{t-L}^R + q_{t-L_0}^E - D_t$, and may be negative since backorder is allowed.  

We now formalize the family of admissible policies $\Pi$, which will determine the new orders placed.
An admissible policy $\pi$ consists of a sequence of measurable maps $\lbrace f^{\pi}_t, t \geq 1 \rbrace$, where each $f^{\pi}_t$ is a deterministic measurable function with domain $\Bbb{R}^{L + L_0 + 1}$ and range $\Bbb{R}^{+,2}$.  In that case, for a given policy $\pi$, the regular order placed in period $t$ equals $f^{\pi}_{R,t}(q_{t-L}^R,\ldots,q_{t-1}^R,q_{t-L_0}^E,\ldots,q_{t-1}^E,I_t)$; while the express order placed in period $t$ equals $f^{\pi}_{E,t}(q_{t-L}^R,\ldots,q_{t-1}^R,q_{t-L_0}^E,\ldots,q_{t-1}^E,I_t)$, and $\Pi$ denotes the family of all such admissible policies $\pi$.

Let $G(y)$ be the sum of the holding and backorder costs when the inventory level equals $y$ in the end of a time period, i.e. $G(y)\stackrel{\Delta}{=}\ hy^+  + by^-$,
where $x^+\stackrel{\Delta}{=} \max(x, 0)$, $x^-\stackrel{\Delta}{=} \max(-x, 0)$.  Here we note that $G$ is convex and Lipschitz, and for $x,y \in \bbr$,
\begin{equation}\label{Lip1}
|G(x) - G(y)| \leq \max(b,h) |x - y|\ \ \ ,\ \ \ |G(x)| \geq \min(b,h) |x|.
\end{equation}

For $t \geq L_0 + 1$, let $C_t$ be the sum of the holding and backorder costs incurred in time period $t$, plus the ordering cost incurred for orders placed in period $t - L_0$, i.e. $C_t \stackrel{\Delta}{=}\ c_R q_{t-L_0}^R +   c_E q_{t-L_0}^E + G(I_t + q_{t-L}^R + q_{t-L_0}^E - D_t).$  We note that charging in period $t$ for orders placed in period $t - L_0$ is a standard ``accounting trick" in the inventory literature to simplify various notations (cf. \citet{Zipkina}), and for the problems considered without loss of generality (w.l.o.g.).  To denote the dependence of the cost on the policy $\pi$, we use the notation $C^{\pi}_t$.  Let $C(\pi)$ denote the long-run average cost incurred by a policy $\pi$, i.e. $C(\pi) \stackrel{\Delta}{=}\  \limsup_{T\to \infty} \frac{\sum_{t = L_0 + 1}^{T} \bbe\left[C^{\pi}_t\right]}{T},$ where we again note that starting the relevant sum at $t  = L_0 + 1$ (as opposed to $t = 1$) is w.l.o.g. for the problems considered.  The value of the corresponding long-run average cost dual-sourcing inventory optimization problem is denoted by $\textrm{OPT}(L) \stackrel{\Delta}{=}\  \inf_{\pi\in \Pi} C(\pi).$

Before proceeding, it will be useful to apply certain well-known reductions to the problem at hand, where we note that similar reductions are known to hold for many classical inventory problems with backlogging (cf. \citet{Karlin58,Scarf60}).  First, as stated in cf. \citet{Sheopuri10}, for the long-run average cost problems which will be the focus of our analysis, any problem with $c_R > 0$ can be transformed into an equivalent problem with $c_R = 0$.  As such we assume throughout that $c_R = 0$.  Let us define the so-called \emph{expedited inventory position} at time $t \geq 1$ as $\hat{I}_t \stackrel{\Delta}{=} I_t + \sum_{k = t - L_0}^{t-1} q_k^E + \sum_{k = t - L}^{t - L + L_0} q_k^R$, which corresponds to the net inventory at the start of period $t$ plus all orders to be received in periods $t,\ldots,t + L_0$ (which were placed before period $t$), and the \emph{truncated regular pipeline} at time $t$ as the $(L - L_0 - 1)$-dimensional vector $\mathcal{R}^t \stackrel{\Delta}{=} (q_{t - L + L_0 + 1}^R,\ldots,q_{t-1}^R)$, with $\mathcal{R}^t_k = q_{t - L + L_0 + k}^R, k = 1,\ldots,L - L_0 - 1$.  Let $\hat{\Pi}$ denote those policies belonging to $\pi$ with the additional restriction that the new orders $q_t^R, q_t^E$ are measurable functions of only $\hat{I}_t, \mathcal{R}^t$.  More formally, $\pi \in \hat{\Pi}$ if there exists a sequence of measurable maps $\lbrace \hat{f}^{\pi}_t, t \geq 1 \rbrace$, where each $\hat{f}^{\pi}_t$ is a deterministic measurable function with domain $\Bbb{R}^{L - L_0}$ and range $\Bbb{R}^{+,2}$, s.t. the regular order placed in period $t$ equals $\hat{f}^{\pi}_{R,t}(\mathcal{R}^t,\hat{I}_t)$ and the express order placed in period $t$ equals $\hat{f}^{\pi}_{E,t}(\mathcal{R}^t,\hat{I}_t)$.

Note that $\hat{I}_1 = - \sum_{i=1}^{\hat{G}} D'_{-i}$, and $\mathcal{R}^1 = \mathbf{0}$.  Also, for any policy $\pi \in \hat{\Pi}$ and $t \geq 1$, it holds that $\hat{I}_{t+1} = \hat{I}_t + q_t^E + \mathcal{R}^t_1 - D_t$, $\mathcal{R}^{t+1}_k = \mathcal{R}^t_{k+1}$ for $k \in [1, L - L_0 - 2]$, and $\mathcal{R}^{t+1}_{L - L_0 - 1} = q_t^R$.  Furthermore, for all $t \geq L_0 + 1$, $C^{\pi}_t = G(\hat{I}_{t - L_0} + q_{t - L_0}^E - \sum_{i = t - L_0}^t D_i) + c q_{t - L_0}^E$.  Then the following is proven in \citet{Sheopuri10}.

\begin{lemma}[\citet{Sheopuri10}\ Lemma 2.1]
$\inf_{\pi\in \Pi} C(\pi) =  \inf_{\pi\in \hat{\Pi}} C(\pi)$, i.e. one may w.l.o.g. restrict oneself to policies belonging to $\hat{\Pi}$.
\end{lemma}
For the remainder of the paper, we thus consider the relevant optimization only over policies belonging to $\hat{\Pi}$, i.e.
\begin{equation}\label{theprob1}
\textrm{OPT}(L) = \inf_{\pi\in \hat{\Pi}} C(\pi).
\end{equation}

For a given policy $\pi \in \hat{\Pi}$, let $\mathcal{R}^{\pi,t} (\hat{I}^{\pi}_t)$ denote a r.v. distributed as the truncated regular pipeline (expedited inventory position) at the start of period $t$ under policy $\pi$.  Similarly, let $q_t^{\pi,E} (q_t^{\pi,R})$ denote the expedited (regular) order placed in period $t$, and suppose that all these r.v.s are constructed on a common probability space, and have the appropriate joint distribution induced by the operation of $\pi$ over time.

\subsection{TBS policy}\label{ds-sec-TBS}
In this section, we formally introduce the family of TBS policies, and characterize the ``best" TBS policy.  A TBS policy $\pi_{r,S}$ with parameters ($r,S$) is defined (cf. JSS) as the policy that places a constant order $r$ from R in every period, and follows an order-up-to rule from E which in each period raises the expedited inventory position to $S$ (if it is below $S$), and otherwise orders nothing.  More formally, under this policy $q^R_t = r$, and $q^E_t = \max(0, S - \hat{I}_t)$, for all $t$.

Let  $I^r_{\infty} \stackrel{\Delta}{=} \sup_{j \geq 0} \left(j r - \sum_{i=1}^j D_{i}\right)$.  In that case, it follows from the results of JSS that
\begin{equation}\label{TBSformula}
C(\pi_{r,S}) =\ c(\bbe[D] - r)+ \bbe\left[G\left( I^r_{\infty} + S - \sum_{i=1}^{L_0 + 1} D'_i \right)\right].
\end{equation}

Note that for each $r$, the minimization problem $\inf_{S \in \bbr} C(\pi_{r,S})$ is equivalent to a standard one-period newsvendor problem.  Furthermore, defining $F^{\infty}(r) \stackrel{\Delta}{=} \inf_{S \in \bbr} C(\pi_{r,S})$, it is proven in JSS that $F^{\infty}(r)$ is convex in $r$ on $\left(-\infty, \bbe[D]\right)$.  Combining the above with standard results for single-server queues (cf. \citet{Asmussen03}) and (\ref{Lip1}), we conclude that there exists at least one pair ($r^*, S^*$) s.t. $r^* \in \argmin_{0\leq r\leq \E[D]} F^{\infty}(r)$ and $S^* \in \argmin_{S \in \bbr} C(\pi_{r^*,S})$; that this pair defines the TBS policy with least long-run-average cost; and that this pair can be computed efficiently by solving a convex program which is independent of the larger lead time $L$.

\subsection{Main result}\label{ds-sec-main}
\subsection{Additional definitions and notations}
Before stating our main result, we will need several additional definitions and notations to describe various relevant quantities which will appear in our bounds on the optimality gap.  For $\theta \geq 0$ and $\epsilon \in \big(0,\E[D]\big]$, let us define 
$$\phi_{\epsilon}(\theta) \stackrel{\Delta}{=} \exp\big(\theta (\E[D] - \epsilon)\big) \bbe[\exp(- \theta D)]\ \ \ ,\ \ \ \gamma_{\epsilon} \stackrel{\Delta}{=} \inf_{\theta \geq 0} \phi_{\epsilon}(\theta),$$
and $\vartheta_{\epsilon} \in \argmin_{\theta \geq 0} \phi_{\epsilon}(\theta)$ denote the supremum of the set of minimizers of $\phi_{\epsilon}(\theta)$, 
where we define $\vartheta_{\epsilon}$ to equal $\infty$ if the above infimum is not actually attained.  Note that $\phi_{\epsilon}(\theta)$ is a continuous and convex function of $\theta$ on $(0,\infty)$, and right-continuous function of $\theta$ at $0$.    
In addition, it follows from \cite{Folland99} Theorem 2.27 that $\phi_{\epsilon}(\theta)$ is right-differentiable at zero, with derivative equal to $- \epsilon$.  We conclude from the definition of derivative and a straightforward contradiction argument that $\vartheta_{\epsilon} > 0$ and $\gamma_{\epsilon} \in [0,1)$.  Let $g \stackrel{\Delta}{=} \inf_{x \in \Bbb{R}} \bbe\left[G\left( x - \sum_{i=1}^{L_0 + 1} D'_i \right)\right] > 0$, and $U \stackrel{\Delta}{=} C(\pi_{0,0}) = c \E[D] + \E[G(- \sum_{i=1}^{L_0 + 1} D'_i)]$, in which case it is easily verified that $g \leq \textrm{OPT}(L) \leq U$ for all $L > L_0 + 1$.  We also make the following additional definitions:
$$p_0 \stackrel{\Delta}{=} \p(D < \E[D]) \in (0,1)\ \ \ ,\ \ \ \hat{p}_0 \stackrel{\Delta}{=} \big(\frac{1}{2} p_0 (1 - p_0)\big)^{\frac{1}{2}} \in (0,1)\ \ \ ,\ \ \ Q_0 \stackrel{\Delta}{=} \inf\lbrace x \in \Bbb{R}^+: \p(D \leq x) \geq \frac{1}{2} p_0 \rbrace \in [0,\E[D]),$$
$$\eta_0 \stackrel{\Delta}{=} \inf_{z \in \Bbb{R}} \E[|z - D|] > 0\ \ \ ,\ \ \ c_0 \stackrel{\Delta}{=} \frac{1}{240} \min(b,h) \hat{p}_0 \eta_0\ \ \ ,\ \ \ U_0 \stackrel{\Delta}{=} 64 (L_0 + 1) \frac{\max^2(b,h)}{\min(b,h)} \E[D],$$
$$\epsilon_0 \stackrel{\Delta}{=} \min\bigg(\E[D] - Q_0 , \frac{1}{4}(\eta_0 \hat{p}_0)^2 , 1 - 2^{-\frac{\hat{p}^2_0}{400}} , \frac{1}{625} c^2_0 \big( U_0 2^{L_0} + \eta_0 + U + 1\big)^{-2} \bigg) \in (0,1 - 2^{-\frac{1}{400}}) \subset (0 , .002),$$
$$Y_0 \stackrel{\Delta}{=} 25 g^{-2}  \big( U_0 2^{L_0} +  \max(b,h) \gamma_{\epsilon_0} \vartheta^{-1}_{\epsilon_0}(1 - \gamma_{\epsilon_0})^{-2} \big)^2 + L_0 + 1.$$
Our main result proves that the best TBS policy is asymptotically optimal as $L\to \infty$, and provides explicit bounds on the optimality gap.
\begin{theorem}\label{ds-thm-main1}
For all $L_0 \geq 0$, $\epsilon \in (0,1)$, and $L > \epsilon_0^{-2} + Y_0 \epsilon^{-2}$, it holds that $\frac{C(\pi_{r^*,S^*})}{\textrm{OPT}(L)} < 1 + \epsilon$.
\end{theorem}
\begin{corollary}\label{maincorr11}
$\lim_{L\to \infty} \frac{C(\pi_{r^*,S^*})}{\textrm{OPT}(L)} = 1.$
\end{corollary}
\section{Proof of Theorem \ref{ds-thm-main1}}\label{ds-sec-proof}
\subsection{Lower bound for the optimal cost}\label{ds-sec-lower}
In this section, we prove a lower bound for $\textrm{OPT}(L)$ by extending the steady-state/convexity approach of \citet{Xin14a} to the dual-sourcing setting.  We note that here our lower bound will involve a non-trivial optimization over measurable functions,
in contrast to the bounds used in \citet{Xin14a} which were of a static nature.  As in \cite{Xin14a}, we will proceed by relating the ``long-run behavior" of ``an optimal policy" to a certain TBS policy.  At a high level, we will combine convexity and the conditional Jensen's inequality with the fact that the r.v.s corresponding to (appropriately defined stationary versions of) the different components of the truncated regular pipeline vector (under the optimal policy) have the same mean, which will (approximately) coincide with the constant order from R in our TBS policy.  Furthermore, when we apply the conditional Jensen's inequality to certain terms corresponding to (appropriately defined stationary versions of) the expedited orders under the same optimal policy, the resulting terms will be suitably measurable functions of past demands, which will (approximately) coincide with the amount of inventory ordered from E in our TBS policy.  
\subsubsection{Connecting to a stationary problem.}
As in \cite{Xin14a}, our program immediately encounters a technical problem.  Namely, the natural way to analyze the ``long-run behavior" of an optimal policy is through the steady-state distribution of the Markov chain induced by this policy.  However, it is not obvious that this steady-state exists.  Actually, it is not even obvious that there exists a stationary optimal policy (so that the dynamics even define a Markov chain), nor even that there even exists an optimal policy at all (as opposed to it only being approached).  Although such questions have been rigorously analyzed for simpler inventory models in \cite{Huh11}, such questions have not been rigorously answered for the setting of more complicated dual-sourcing models.  We note that although in \cite{Sheopuri10} it is stated in passing that many of the same results should extend to the dual-sourcing setting, no proofs are provided, and the explicit assumptions needed for such a transference are not clarified.  A similarly terse exposition on related questions is provided in \cite{Hua14}.  Furthermore, in none of these works is the question of existence of and convergence to relevant stationary measures discussed.  To overcome this, as in \cite{Xin14a}, we first observe that we will not actually need a random vector which is truly the steady-state of the aforementioned Markov chain (which in principle may not exist), but only need to demonstrate the existence of a random vector which has several properties that we would want such a steady-state (if it existed) to have.  We now show the existence of such a random vector.  We note that although closely related questions have been studied in the MDP literature (cf. \cite{Arapostathis93}), and perturbative approaches similar to the approach we take in our own proof are in general well-known (cf. \cite{Filar}), to the best of our knowledge the desired result does not follow directly from any results appearing in the literature.  As such, we include a proof for completeness in the technical appendix Section\ \ref{techsec}.  We note that here the relevant analysis is considerably more challenging than that given in \cite{Xin14a}, due to the fact that in the dual-sourcing setting the inventory level is unbounded from below, and the associated ordering levels are not known to be uniformly bounded (in contrast to the setting considered in \cite{Xin14a} for which such bounds were already proven in \cite{Zipkina}).  Furthermore, although several bounds exist in the dual-sourcing literature relating order levels under an optimal policy to the inventory level at the time of ordering (cf. \cite{Sheopuri10,Hua14}), it seems that due to the inventory being unbounded below none of those bounds are suitable for our purposes.  It is also worth noting that our approach is able to side-step many of the complexities and additional assumptions (e.g. finite second moment or bounded support) often required when analyzing inventory models which are unbounded from below. 
\\\indent We defer all relevant proofs to the technical appendix Section\ \ref{techsec}.  For two r.v.s $X,Y$, let $X \sim Y$ denote equivalence in distribution.  Before stating our result, for the sake of building intuition, we first describe what the various r.v.s appearing in our result would correspond to ``if we were to assume" (which we do not, i.e. it is not an assumption of our main results) that there exists an optimal policy which is stationary, and whose corresponding Markov chain converges to a steady-state distribution, i.e. the truncated regular pipeline and expedited inventory position converge in distribution under the operation of this optimal stationary policy.  In that case, our theorem contains an $(L - L_0 - 1)$-dimensional random vector $\mathbf{\chi}^{*, L}$, an $(L - L_0)$-dimensional random vector $\mathbf{q}^{*,L}$, and a r.v. $\mathcal{I}^{*, L}$, which may be interpreted as follows.  Suppose one has been operating under this stationary optimal policy for a long time, say up to some very large time $T$, at which time the system is essentially in steady-state (again we note that this discussion is purely for the sake of building intuition, and our main results do not actually assume this).  Then $\mathbf{\chi}^{*, L}$ corresponds to the steady-state truncated regular pipeline vector under this optimal policy (at time $T$), i.e. $\mathbf{\chi}^{*, L}_i$ is the regular order which enters the expedited inventory position in period $T + i - 1$.  $\mathbf{q}^{*,L}$ corresponds to the steady-state vector of expedited orders to be placed over the next $L - L_0$ periods under this optimal policy, i.e. $\mathbf{q}^{*, L}_i$ is the expedited order which enters the expedited inventory position in period $T + i - 1$.  Finally, $\mathcal{I}^{*, L}$ corresponds to the steady-state expedited inventory position under this optimal policy (at time $T$).

\begin{theorem}\label{proxychain0}
For all $L_0 \geq 0$ and $L > L_0 + 1$, one may construct  an $L - L_0 - 1$-dimensional random vector $\mathbf{\chi}^{*, L}$, an $L - L_0$-dimensional random vector $\mathbf{q}^{*, L}$, and a random variable $\mathcal{I}^{*, L}$, as well as $\lbrace D_i, i \geq 1 \rbrace$, on a common probability space s.t. the following are true.
\begin{enumerate}[(i)]
\item \label{enum0} W.p.1 $(\mathbf{\chi}^{*, L}, \mathbf{q}^{*, L})$ is non-negative.  Also, $(\mathbf{\chi}^{*, L},\mathcal{I}^{*, L})$ is independent of $\lbrace D_i, i \geq 1 \rbrace$, and $q^{*, L}_i$ is independent of $\lbrace D_j, j \geq i \rbrace$ for $i \in [1, L - L_0]$.
\item \label{enum0b} $\mathbf{\chi}^{*, L}_i \sim \mathbf{\chi}^{*, L}_1$ for $i \in [1 , L - L_0 - 1]$, and $q^{*, L}_i \sim q^{*, L}_1$ for $i \in [1 , L - L_0]$.
\item \label{lseq-300} For all $k \in [1,L - L_0]$,
$$\mathcal{I}^{*, L} + \sum_{i=1}^{k-1} (q^{*, L}_i + \chi^{*, L}_i - D_i) + q^{*, L}_k - \sum_{i = k}^{k + L_0} D_i
\sim \mathcal{I}^{*, L} + q^{*, L}_1 - \sum_{i=1}^{L_0 + 1} D_i.$$
\item \label{enum0mean1}$(\mathbf{\chi}^{*, L}, \mathbf{q}^{*, L},\mathcal{I}^{*, L})$ has finite mean.
\item \label{enum0mean2} $\E[\mathbf{\chi}^{*, L}_1] + \E[q^{*, L}_1] = \E[D]$.
\item \label{enum0mean3} $$\textrm{OPT}(L) \geq c\left(\bbe[D] - \E[\chi^{*, L}_1]\right) + \bbe \left[G\left( \mathcal{I}^{*, L} + q^{*, L}_1 - \sum_{i=1}^{L_0 + 1} D_i\right)\right].$$
\end{enumerate}
\end{theorem}

\subsubsection{Vanishing discount factor approach.}
Although Theorem\ \ref{proxychain0}.(\ref{enum0mean3}) relates $\textrm{OPT}(L)$ to a certain expectation, this expectation (as written) is not immediately amenable to analysis.  To remedy this, we introduce a discount factor $\alpha$ to implement the so-called ``vanishing discount factor" approach to analyzing infinite-horizon MDP (cf. \citet{Huh11}), which will allow for a simpler analysis when we pass to the limit as $L \rightarrow \infty$.  Indeed, this discount factor will help us to analyze the lower bound which arises when we apply the conditional Jensen's inequality, as this lower bound will itself involve the solution to a non-trivial multi-stage dynamic optimization problem.  We note that the lower bound which arose when related techniques were applied to single-sourcing systems with lost sales in \citet{Xin14a} only involved a static optimization problem, and thus no such discount factor was introduced.  In particular, Theorem\ \ref{proxychain0} immediately implies the following corollary.  Let $r_L \stackrel{\Delta}{=} \E[\chi^{*, L}_1]$.

\begin{corollary}\label{proxy00}
For all $L_0 \geq 0, L > L_0 + 1$, and $\alpha \in (0,1)$, 
\begin{eqnarray*}
\textrm{OPT}(L) &\geq& c\left(\bbe[D] - r_L\right) + \frac{1-\alpha}{1-\alpha^{L}} \sum_{k=1}^{L} \alpha^{k-1}  \bbe \left[G\left( \mathcal{I}^{*, L} + q^{*, L}_1 - \sum_{i=1}^{L_0 + 1} D_i\right)\right]
\\&\geq& c\left(\bbe[D] - r_L\right) + (1 - \alpha) \sum_{k = 1}^{L - L_0} \alpha^{k - 1}  \bbe \left[G\left(\mathcal{I}^{*, L} + \sum_{i=1}^{k-1} (q^{*, L}_i + \mathbf{\chi}^{* , L}_i - D_i) + q^{*, L}_k - \sum_{i = k}^{k + L_0} D_i\right)\right].
\end{eqnarray*}
\end{corollary}

\subsubsection{Applying the conditional Jensen's inequality and relating to a single-source inventory model.}\label{3pp}
We now apply the conditional Jensen's inequality to Corollary\ \ref{proxy00}, which will allow us to lower-bound $\textrm{OPT}(L)$ by the optimal value of a certain finite-horizon single-source inventory model with backlogged demand.  We will then relate this finite-horizon problem to an associated infinite-horizon problem, which has an optimal stationary policy.  Furthermore, we will connect the behavior of such an optimal stationary policy to the performance of an associated TBS policy, ultimately allowing us to prove our main results.
In particular, it follows from Theorem\ \ref{proxychain0} and the independence structure of the relevant r.v.s that for $k \in [1, L - L_0]$,
$$
\E\bigg[ \mathcal{I}^{*, L} + \sum_{i=1}^{k-1} (q^{*, L}_i + \mathbf{\chi}^{* , L}_i - D_i) + q^{*, L}_k - \sum_{i = k}^{k + L_0} D_i
\bigg| D_{[k + L_0]}\bigg]
$$
equals
$$
\E[\mathcal{I}^{*, L}] + \sum_{i=1}^{k-1} (\E[q^{*, L}_i | D_{[i-1]}] + r_L - D_i) + \E[q^{*, L}_k | D_{[k-1]}] - \sum_{i = k}^{k + L_0} D_i.$$
Further combining with Corollary\ \ref{proxy00}, the convexity of $G$, and Jensen's inequality for conditional expectations, we obtain the following result.
\begin{proposition}\label{OBSuse1}
For any $\alpha\in (0,1)$ and $L> L_0 + 1$, $\textrm{OPT}(L) -  c \left(\bbe[D] - r_L \right)$ is at least
\begin{equation}\label{refer1bb}
\begin{aligned}
& (1 - \alpha) \sum_{k = 1}^{L - L_0} \alpha^{k-1} \bbe\Bigg[G\bigg(
\E[\mathcal{I}^{*, L}] - (L_0 + 1) r_L + \sum_{i=1}^{k-1} \big(\E[q^{*, L}_i | D_{[i-1]}] - (D_i - r_L)\big) \\
&\ \ \ \ \ \ \ \ \ \ \ \ \ \ \ \ \ \ \ \ \ \ \ \ \ \ \ \ \ \ \ \ \ \ \ \ \ \ \ \ \ \ \ \ \ \ \ \ \ \ \ \      +\E[q^{*, L}_k | D_{[k-1]}] - \sum_{i = k}^{k + L_0} (D_i - r_L)
\bigg)\Bigg].
\end{aligned}
\end{equation}
\end{proposition}

Note that (\ref{refer1bb}) is the discounted cost incurred (during periods $L_0 + 1,\ldots,L$) by the policy ordering $\E[q^{*, L}_i | D_{[i-1]}]$ in period $i$, of a single-sourcing $L$-period backlog inventory problem with unit holding cost $h$, backorder cost $b$, zero ordering cost, discount factor $\alpha$, i.i.d. demand distributed as $D - r_L$ (which we note can be positive \emph{or negative}), lead time $L_0$, and initial inventory position (initial net inventory plus all entries of the initial pipeline vector) $\E[\mathcal{I}^{*, L}] - (L_0 + 1) r_L$ (cf. \citet{Karlin58}), multiplied by $(1 - \alpha)$.  Such models, and their optimal policies, have been studied in-depth in the literature (cf. \citet{Karlin58, Zipkin00, Fleischmann03}), and are well-understood (especially for the case of non-negative demand, cf. \citet{Zipkin00}).  Let $\overline{\Pi}$ denote the family of all feasible non-anticipative policies for the aforementioned inventory problem (as it is typically defined, cf. \citet{Zipkin00}).  For $\pi \in \overline{\Pi}$, initial inventory position $x \in \bbr$, $r \in \bbr$, and $i \geq 1$, let $C^{\pi}_i(r,x)$ denote the cost incurred by policy $\pi$ in the aforementioned inventory problem in period $i + L_0$, if the demand in each period is i.i.d. distributed as $D - r$ (with the leadtime $L_0$ and costs $b,h$ as above).  For $x \in \bbr, r \in \bbr, \alpha \in (0,1), n \geq 1$, let us define
\begin{equation}\label{ds-finite}
V_{\alpha}^{n}(r,x) \stackrel{\Delta}{=}\ \inf_{\pi\in \overline{\Pi}} \bbe\left[\sum_{i = 1}^n \alpha^{i-1} C^{\pi}_i(r,x)\right];
\end{equation}
and
\begin{equation}\label{ds-infinite}
V_{\alpha}^{\infty}(r,x) \stackrel{\Delta}{=}\ \inf_{\pi\in \overline{\Pi}} \bbe\left[\sum_{ i = 1}^{\infty} \alpha^{i-1} C^{\pi}_i(r,x)\right].
\end{equation}
As a notational convenience, we define $V_{\alpha}^0(r,x) = 0$, $V_{\alpha}^n(r,-\infty) \stackrel{\Delta}{=} \inf_{x \in \bbr} V_{\alpha}^n(r,x), V_{\alpha}^{\infty}(r,-\infty) \stackrel{\Delta}{=} \inf_{x \in \bbr} V_{\alpha}^{\infty}(r,x)$.  Then combining the above, we derive the following lower bound for $\textrm{OPT}(L)$.

\begin{lemma}\label{reduce0}
For all $L_0 \geq 0, L > L_0 + 1$, and $\alpha \in (0,1)$, 
\begin{equation}\label{inter1}
\textrm{OPT}(L) \geq c (\bbe[D] - r_L) + (1 - \alpha) V_{\alpha}^{L - L_0}(r_L, -\infty).
\end{equation}
\end{lemma}

\subsubsection{Overview of remainder of the proof of our main results.}

The remainder of the proof involves a careful analysis of the right-hand-side (r.h.s.) of (\ref{inter1}) as $L \rightarrow \infty$, and we now sketch an outline of our approach.  First, we will prove that if $r_L$ is bounded away from $\E[D]$, then 
$V_{\alpha}^{\infty}(r_L, -\infty) - V_{\alpha}^{L - L_0}(r_L, -\infty)$ can be suitably bounded by a function of $L$ which converges to 0 as $L \rightarrow \infty$.  We then observe that the infinite-horizon problem associated with $V_{\alpha}^{\infty}(r_L, -\infty)$ has an optimal policy which is stationary, Markov, and of order-up-to type.  Furthermore, the stochastic process induced  by this optimal policy will be equivalent to that induced by a corresponding TBS policy, but possibly initialized \emph{not according to the stationary distribution of the associated inventory process}.  Then we prove that $r_L$ is indeed bounded away from $\E[D]$, since otherwise we can use the theory of random walks to derive a contradiction (as $\textrm{OPT}(L)$ would be strictly greater than $U$).  Finally, we combine these facts to bound various error terms under a suitable choice of $\alpha$ (which converges to 1 as $L \rightarrow \infty$), including a term resulting from the difference between the performance of the same TBS policy under different initializations, to prove our main result Theorem\ \ref{ds-thm-main1}.

\subsection{Proof of Theorem\ \ref{ds-thm-main1}}
We now complete the proof of Theorem\ \ref{ds-thm-main1} by formalizing the argument sketched at the end of Section\ \ref{ds-sec-lower}.  Such arguments are standard in the literature on MDP and infinite-horizon inventory control problems (cf. \citet{Iglehart63, Sennott89, Schal93, Fleischmann03, Feinberg11, Huh11}).  We note that the somewhat non-standard aspect here is that the demand in each period is distributed as $D - r_L$, and thus \emph{may be negative}.  As such, the original arguments typically used to analyze the relevant quantities and prove related interchange-of-limits results (cf. \citet{Iglehart63}) do not directly apply.  The possibility of negative demand also makes the verification of the conditions of general theorems which validate such bounds and interchange-of-limits (cf. \citet{Sennott89,Schal93}) somewhat involved, even when these theorems are customized to the inventory setting (cf. \citet{Parker04, Huh11}).  We note that the verification of closely related results have arisen recently in the context of analyzing inventory systems with returns, which reduce to standard inventory systems where demand can be positive or negative (cf. \citet{Fleischmann03}).  However, those results (which verify the technical conditions of \citet{Sennott89}) do not seem to extend immediately to our case, and further seem to require that the demand and ordering quantities take integer values.  In light of the above, and for the sake of clarity and completeness, we now provide a self-contained proof of all necessary bounds, which (combined with Lemma\ \ref{reduce0}) will complete the proof of our main result Theorem\ \ref{ds-thm-main1}.  We defer most proofs to the technical appendix Section\ \ref{techsec}.

We begin by stating some well-known properties of $V_{\alpha}^{n}(r,x)$ and $V_{\alpha}^{\infty}(r,x)$, which follow from the results of JSS, \citet{Karlin58} and \citet{Scarf60}.  We note that although in some cases the proofs there are only explicitly given for the case of non-negative demand, as noted in \citet{Heyman84} and \citet{Fleischmann03}, the arguments carry over to the general case (in which demand may be negative) with only trivial modification.

\begin{lemma}[JSS, \citet{Scarf60}]\label{use1}
For all $\alpha \in (0,1), r, x \in \bbr$, and $n \geq 1$,
$$
V_{\alpha}^{n}(r,x) = \inf_{y \geq x} \left( \bbe\left[G\big(y - \sum_{k=n}^{L_0 + n} (D_k - r) \big)\right] + \alpha \bbe\left[V_{\alpha}^{n-1}\big(r, y - (D_{L_0 + n} - r)\big)\right] \right).$$
Furthermore, $V_{\alpha}^n(r,x)$ is: a convex (and thus also continuous) function of $x$ on $\bbr$ for each fixed $n,r$; a continuous function of $r$ on $\bbr$ for each fixed $n,x$; an increasing function of $x$ on $\bbr$ for each fixed $n,r$; and an increasing function of $n$ on $Z^+$ for each fixed $x,r$.  In addition, the infinite-horizon problem stated in the r.h.s. of (\ref{ds-infinite}) admits an optimal stationary Markov policy.
\end{lemma}

Next, we bound $V_{\alpha}^{\infty}(r,x) - V_{\alpha}^n(r,x)$, and combine our bounds with Lemma\ \ref{use1} to derive some useful properties of $V_{\alpha}^{\infty}(r,x)$ and the associated optimization problem.  We defer all proofs to the technical appendix Section\ \ref{techsec}.  Let $\overline{S}_{\alpha}(r) \stackrel{\Delta}{=} 4 (L_0 + 1) \frac{\max(b,h)}{\min(b,h)}(|r| + \E[D])(1 - \alpha)^{-2}$.

\begin{lemma}\label{limequal2}
For $\alpha \in (0,1)$, $r, x \in \bbr$, and $n \geq 1$,
\begin{equation}\label{limequal2a}
0 \leq  V_{\alpha}^{\infty}(r,x) - V_{\alpha}^n(r,x) \leq \max(b,h) \big(\overline{S}_{\alpha}(r) + |x| + |r| + \E[D]\big) (1 + L_0 + n) (1 - \alpha)^{-2} \alpha^n,
\end{equation}
and $V_{\alpha}^{\infty}(r,x) = \lim_{n \rightarrow \infty} V_{\alpha}^n(r,x)$.  Furthermore, for $\alpha \in (0,1)$ and $r \in \bbr$, $V_{\alpha}^{\infty}(r,x)$ is a finite-valued, convex, and non-decreasing function of $x$ on $\bbr$.  Letting $S_{\alpha}^{\infty}(r)$ denote the supremum of the set of minimizers (in $x$) of $V_{\alpha}^{\infty}(r,x)$, it holds that $|S_{\alpha}^{\infty}(r)| \leq \overline{S}_{\alpha}(r)$, and the infinite-horizon problem stated in the r.h.s. of (\ref{ds-infinite}) admits an optimal stationary base-stock policy, with order-up-to level $S_{\alpha}^{\infty}(r)$.  In addition, for $L_0 \geq 0$, $L > L_0 + 1$, and $\alpha \in (0,1)$,
\begin{equation}\label{limequal2b}
\textrm{OPT}(L) \geq c(\bbe[D] - r_L) + (1 - \alpha) V_{\alpha}^{\infty}\big(r_L, S_{\alpha}^{\infty}(r_L) \big) -  U_0 (1 - \alpha)^{-3} L \alpha^{L - L_0}.
\end{equation}
\end{lemma}

We now formally define the Markov process representing the inventory position process under such an optimal stationary base-stock policy, initialized in state $S_{\alpha}^{\infty}(r_L)$.  Let $S_{\alpha,L} \stackrel{\Delta}{=} S_{\alpha}^{\infty}(r_L)$.  For $r \in [0,\E[D]]$ and $y \in \Bbb{R}$, let $\lbrace X^{r,y}_k, k \geq 1 \rbrace$ denote the following Markov process.  $X^{r,y}_1$ equals $y$.  For all $k \geq 1$, $X^{r,y}_{k+1} = \max\big( X^{r,y}_k + r - D_k, y \big)$.  Let $W^r_k \stackrel{\Delta}{=} \sum_{j=1}^k (r - D_j)$,
$Z^r_k \stackrel{\Delta}{=} \max_{i \in [0,k - 1]} W^r_i$, $Z^r_{\infty} \stackrel{\Delta}{=}  \sup_{i \geq 0} W^r_i$, $M^r_k \stackrel{\Delta}{=} \E[Z^r_k]$, $M^r_{\infty} \stackrel{\Delta}{=} \E[Z^r_{\infty}]$.  It follows from the well-known analysis of the single-server queue using Lindley's recursion (cf. \citet{Asmussen03}) that $X^{r,y}_k \sim y + Z^r_k$; and $X^{r,y}_{\infty} \stackrel{\Delta}{=} \lim_{k \rightarrow \infty} X^{r,y}_{k}$ (in the sense of weak convergence) is a well-defined r.v. distributed as $y + Z^r_{\infty}$.  

Combining these definitions with Lemmas\ \ref{use1} and\ \ref{limequal2}, we conclude the following.

\begin{corollary}\label{wrapit0}
For $L_0 \geq 0, L > L_0 + 1$, and $\alpha \in (0,1)$,
$$
\textrm{OPT}(L) \geq c(\bbe[D] - r_L) + (1 - \alpha) \sum_{k=1}^{\infty} \alpha^{k-1} \E\big[G\big(S_{\alpha,L} + Z^{r_L}_k - \sum_{i=1}^{L_0 + 1} (D'_i - r_L)\big)\big] -  U_0 (1 - \alpha)^{-3} L \alpha^{L - L_0}.
$$
\end{corollary}

We now briefly review some useful properties of $Z^r_k$, which we will use to complete the proof of our main results.  These properties follow by combining generally well-known results for generating functions, large deviations, single-server queues, and recurrent random walks
(cf. \citet{Spitzer56, Kingman62a, Folland99, Asmussen03, Xin14a}), and we omit the details.  

\begin{lemma}\label{zreview}
For all $r > 0$, $\lbrace M^r_k, k \geq 1 \rbrace$ is non-decreasing, $M^r_{\infty} = \lim_{k \rightarrow \infty} M^r_k$, and for all $i \geq j \geq 1$, $M^r_i - M^r_j = \sum_{k=j}^{i-1} k^{-1} \E[\max(0,W^r_k)].$  
If there exists $\epsilon \in (0, \E[D])$ s.t. $r \leq \E[D] - \epsilon$, then $M^r_{\infty} < \infty$, and $M^r_{\infty} - M^r_{n} \leq \big(\vartheta_{\epsilon}(1 - \gamma_{\epsilon})\big)^{-1} \gamma_{\epsilon}^n$ for all $n \geq 1$.  
\end{lemma}

Finally, we will also need the following corollary (of Lemma\ \ref{zreview}), which shows that $r_L$ is uniformly bounded away from $\E[D]$ in an appropriate sense, and whose proof we again defer to the technical appendix Section\ \ref{techsec}.
\begin{corollary}\label{lbound1}
For all $L > \epsilon_0^{-2} + L_0 + 1$, it holds that $r_L < \E[D] - \epsilon_0$.
\end{corollary}

We now complete the proof of our main results.
\proof{Proof of Theorem\ \ref{ds-thm-main1}}
It follows from (\ref{TBSformula}) that for all $\alpha \in (0,1)$,
\begin{eqnarray*}
C(\pi_{r_L,S_{\alpha,L} + (L_0 + 1) r_L}) &=& c (\E[D] - r_L) + \bbe\left[G\left( S_{\alpha,L} + Z^{r_L}_{\infty} + (L_0 + 1) r_L - \sum_{i=1}^{L_0 + 1} D'_i\right)\right]
\\&=& c (\E[D] - r_L) + (1 - \alpha) \sum_{k=1}^{\infty} \alpha^{k-1} \bbe\left[G\left( S_{\alpha,L} + Z^{r_L}_{\infty} -  \sum_{i=1}^{L_0 + 1} (D'_i - r_L) \right)\right].
\end{eqnarray*}
Combining with Corollaries\ \ref{wrapit0} and\ \ref{lbound1}, Lemma\ \ref{zreview}, (\ref{Lip1}), and the fact that $L > \epsilon_0^{-2} + L_0 + 1$, we conclude that for all $\alpha \in (0,1)$, $C(\pi_{r^*,S^*}) - \textrm{OPT}(L) - U_0 (1 - \alpha)^{-3} L \alpha^{L - L_0}$ is at most
\begin{eqnarray*}
\ &\ &\ (1 - \alpha) \sum_{k=1}^{\infty} \alpha^{k-1}
\left(
\bbe\left[G\left(S_{\alpha,L} + Z^{r_L}_{\infty} - \sum_{i=1}^{L_0 + 1} (D'_i - r_L) \right)\right]
-
\E\left[G\left(S_{\alpha,L} + Z^{r_L}_k - \sum_{i=1}^{L_0 + 1} (D'_i - r_L)\right)\right]
\right) \nonumber
\\&\ &\ \ \ \ \ \ \leq\ \ \ \max(b,h) (1 - \alpha) \sum_{k=1}^{\infty} \alpha^{k-1} \big( \vartheta_{\epsilon_0}(1 - \gamma_{\epsilon_0})\big)^{-1} \gamma_{\epsilon_0}^k \nonumber
\\&\ &\ \ \ \ \ \ =\ \ \ \max(b,h) \gamma_{\epsilon_0} \big(\vartheta_{\epsilon_0}(1 - \gamma_{\epsilon_0})\big)^{-1} \frac{1 - \alpha}{1 - \gamma_{\epsilon_0} \alpha} \nonumber
\\&\ &\ \ \ \ \ \ \leq\ \ \ (1 - \alpha) \max(b,h) \gamma_{\epsilon_0} \vartheta^{-1}_{\epsilon_0}(1 - \gamma_{\epsilon_0})^{-2}. \nonumber
\end{eqnarray*}
We conclude that for all $\alpha \in (\frac{1}{2},1)$, $C(\pi_{r^*,S^*}) - \textrm{OPT}(L)$ is at most
$$
U_0 2^{L_0} (1 - \alpha)^{-3} L \alpha^L + (1 - \alpha) \max(b,h) \gamma_{\epsilon_0} \vartheta^{-1}_{\epsilon_0}(1 - \gamma_{\epsilon_0})^{-2}.$$
As $L > \epsilon_0^{-2} + L_0 + 1$ implies $L > 100$, which itself may be shown to imply that $5 \frac{\log(L)}{L} < \frac{1}{2}$, we may set $\alpha = 1 - 5 \frac{\log(L)}{L}$.  Then applying the fact that $1 - \alpha \leq \exp(- \alpha)$, we conclude that
$$C(\pi_{r^*,S^*}) - \textrm{OPT}(L) \leq 5 \frac{\log(L)}{L} \big( U_0 2^{L_0} +  \max(b,h) \gamma_{\epsilon_0} \vartheta^{-1}_{\epsilon_0}(1 - \gamma_{\epsilon_0})^{-2} \big).$$
As $\frac{\log(L)}{L} < L^{-\frac{1}{2}}$ for all $L \geq 1$, combining with the fact that $\textrm{OPT} \geq g$ and a straightforward calculation completes the proof.  $\Halmos$.
\endproof

\section{Conclusion}\label{ds-sec-conclusion}
In this paper, we proved that when the lead time of the express source is held fixed, a simple TBS policy is asymptotically optimal for the dual-sourcing inventory problem as the lead time of the regular source grows large.   Our results provide a solid theoretical foundation for several conjectures and numerical experiments appearing previously in the literature regarding the good empirical performance of such policies.   Furthermore, the simple TBS policy performs nearly optimally exactly when standard DP-based methodologies become intractable due to the curse of dimensionality.  In addition, since the ``best" TBS policy can be computed by solving a convex program that does not depend on the lead time of the regular source, and is easy to implement, our results lead directly to very efficient algorithms with asymptotically optimal performance guarantees.  We also explicitly bound the optimality gap of the TBS policy for any fixed lead time (of the regular source), and prove that this decays inverse-polynomially in the lead time of the regular source.  Perhaps most importantly, since many companies are already implementing such TBS policies, our results provide strong theoretical support for the widespread use of TBS policies in practice.

This work leaves many interesting directions for future research.  First, it would be interesting to further investigate the rate of convergence to optimality of TBS policies as the lead time grows large, especially in light of their use in practical settings.  Although we have not optimized the explicit bounds which we have proven on the optimality gap, we suspect that proving significantly stronger (e.g. exponentially decaying) bounds will require the development of new techniques.  For example, when we apply the conditional Jensen's inequality to lower bound the optimal value by a certain single-sourcing problem in Section\ \ref{3pp}, our current approach does not incorporate the fact that $\E[q^{*, L}_i]$ is the same for all $i$, instead only using the fact that $\E[q^{*, L}_i | D_{[i-1]}]$ is a measurable function of $D_{[i-1]}$.  It seems plausible that incorporating this ``stationary expectations" property may be a promising approach here.  Previous bounds from the literature on the rate of convergence of finite horizon inventory optimization problems to their infinite horizon counterparts, e.g. \citet{Hordijk74,Hordijk75}, may also be helpful.

Second, and related to the aforementioned discussion as regards the rate of convergence to optimality of TBS policies, it would be interesting to identify other more sophisticated algorithms which perform better for small-to-moderate lead times, yet remain efficient to implement.  Indeed, it remains an interesting open question to better understand the trade-off between algorithmic run-time and acheivable performance guarantees in this context, i.e. how complex an algorithm is required to ``exploit" the weak correlations which persist even as the lead time grows large.  In the context of dual-sourcing, potential algorithms here include: the so-called dual-sourcing smoothing policies recently studied in \citet{Boute14}; affine policies more generally (cf. \citet{Bertsimas10}), of which dual-sourcing smoothing policies are a special case; the single index and dual index policies discussed earlier; or the dual-balancing policies analyzed in \citet{Levi08}.  It would also be quite interesting to analyze ``hybrid" algorithms, which could e.g. solve a large dynamic program when the lead time is small, and gradually transition to using simpler heuristics as the lead time grows large; or combine different heuristics depending on the specific problem parameters.  In the context of the above conversation on optimality gaps, we do remind the reader that for any fixed regular lead time a TBS policy is not exactly optimal except in some very special cases (cf. JSS), and that our results (and associated insights) should always be applied with care.

On a final note, combined with the results of \citet{Goldberg12} and \citet{Xin14a}, our methodology lays the foundations for a completely new approach to analyzing inventory models with large lead times.  So far, this approach has been successful in yielding key insights and efficient algorithms for two settings previously believed intractable: lost sales models with large lead times, and dual-sourcing models with large lead time gap.  We believe that our techniques have the potential to make similar progress on many other difficult supply chain optimization problems of practical relevance in which there is a lag between when policy decisions are made and when those decisions are implemented.  This includes both more realistic variants of the lost-sales and dual-sourcing models considered so far (e.g. models with distributional dependencies, parameter uncertainty, complex network structure, and more accurate modeling of costs), as well as fundamentally different models (e.g. inventory systems with remanufacturing when the manufactured and remanufactured lead times differ, cf. \citet{Zhou11}; multi-echelon systems with lost sales and positive lead times, cf. \citet{Huh10}; or models with perishable goods).  In closing, we note that our approach can more generally be viewed as a methodology to formalize the notion that when there is a high level of uncertainty and randomness in one's supply chain, even simple policies perform nearly as well as very sophisticated policies, since no algorithm can ``beat the noise".  Exploring this concept from a broader perspective may be fruitful in yielding novel algorithms and insights for a multitude of problems in operations management and operations research.

\bibliographystyle{nonumber}

\section{Technical Appendix}\label{techsec}
\subsection{Proof of Theorem\ \ref{proxychain0}}
\subsubsection{Overview of proof.}
Before providing a formal proof, we first provide an intuitive overview, noting that our proof is similar to several proofs in the literature (cf. \cite{Xin14a} and the references therein).  We proceed by constructing a sequence of random vectors, one for each sufficiently small $\epsilon > 0$, and later take an appropriate weak limit (which will become the vector satisfying the conditions of the theorem).  As in \cite{Xin14a}, given $\epsilon > 0$, we will pick a sufficiently large time $T_{\epsilon}$ s.t. the expected performance of an approximately optimal (possibly non-stationary) policy $\pi^{*,\epsilon}$ up to time $T_{\epsilon}$ is ``close" to $\textrm{OPT}(L)$.  We then further prove the existence of a time $T_{1,\epsilon}$ ``near" $T_{\epsilon}$ s.t. the expedited inventory position and truncated regular pipeline vector (under policy $\pi^{*,\epsilon}$) are ``well-behaved" at time $T_{1,\epsilon}$, which will be necessary for our later arguments, as it will allow us to bound the time needed to ``clear the system" if one orders nothing from that time onwards.  We then construct a ``modified policy" and associated Markov chain, which behaves exactly like the expedited inventory position and truncated regular pipeline vector under $\pi^{*,\epsilon}$ on $[1,T_{1,\epsilon}]$, but after that time forces a sequence of ordering decisions which cause the associated inventory position and pipeline vector to re-enter a state distributed as its initial state, at which time the entire process restarts.  We note that due to the process being unbounded from below, here the special initialization involving $-\sum_{i=1}^{\hat{G}} D'_{-i}$ will prove useful.  This regenerative structure, combined with our careful selection of $T_{1,\epsilon}$, will allow us to apply the theory of regenerative processes to prove the existence of a stationary distribution for the associated Markov chain, which we will prove to satisfy the conditions of an ``approximate" version of Theorem\ \ref{proxychain0} (with the approximation error parametrized by $\epsilon$).  Taking a weak limit (as $\epsilon \downarrow 0$) of the associated sequence of random vectors yields a random vector satisfying the conditions of Theorem\ \ref{proxychain0}, completing the proof.
\\\indent As all results in this subsection will be stated for a fixed $L_0 \geq 0$ and $L > L_0 + 1$, we assume these parameters are fixed and supress any associated notational dependencies.  For $\epsilon > 0$, let $\pi^{*,\epsilon}$ denote some (fixed) policy in $\hat{\Pi}$ s.t. $\textrm{OPT}(L) > C(\pi^{*,\epsilon}) - \frac{\epsilon}{2}$.  Let $U_{\epsilon} \stackrel{\Delta}{=} L^2(\frac{4 U}{\epsilon \min(b,h)} + 2 \E[D])$, and $U_{2,\epsilon} \stackrel{\Delta}{=} \big(L + \frac{U_{\epsilon}}{\E[D]} + 2\big) (h + b + c) U_{\epsilon}$.  It follows from the definition of $\limsup$ that there exists $T_{\epsilon} > 100 \big( U_{2,\epsilon} + (U+1) L \big) \epsilon^{-1}$ s.t. $C(\pi^{*,\epsilon}) > T^{-1} \sum_{t= L_0 + 1}^{T} \bbe\left[C^{\pi^{*,\epsilon}}_t \right]  - \frac{\epsilon}{2}$ for all $T \geq (1 - \epsilon) T_{\epsilon} - L$.  
\subsubsection{Existence of time $T_{1,\epsilon}$, near $T_{\epsilon}$, at which inventory and pipeline are small in expectation.}
We first prove that there must exist a time ``close to" $T_{\epsilon}$ at which the expedited inventory position and truncated regular pipeline vector (under policy $\pi^{*,\epsilon}$) are ``small" (in absolute value) in expectation.  
\begin{claim}\label{goodtime1}
For all $\epsilon \in \big(0, \min(\frac{1}{2}, \frac{U}{2}) \big)$, there exists $T_{1,\epsilon} \in  [(1 - \epsilon) T_{\epsilon} - L, T_{\epsilon}]$ s.t. 
\begin{equation}\label{goodtimea}
C(\pi^{*,\epsilon}) > T_{1,\epsilon}^{-1} \sum_{t = L_0 + 1}^{T_{1,\epsilon}} \bbe\left[C^{\pi^{*,\epsilon}}_t \right]  - \frac{\epsilon}{2};
\end{equation}
for all $k \in [0,L - 1]$, 
\begin{equation}\label{goodtime2}
\E[|\hat{I}^{\pi^{*,\epsilon}}_{T_{1,\epsilon} + k - L_0} + q^{\pi^{*,\epsilon},E}_{T_{1,\epsilon} + k - L_0}|] \leq \frac{2 U L}{\epsilon \min(b,h)} + (L_0 + 1) \E[D];
\end{equation}
and for all $k \in [1, L - L_0 - 1]$,
\begin{equation}\label{goodtime3}
\E[{\mathcal R}^{\pi^{*,\epsilon},T_{1,\epsilon} - L_0}_k] \leq 
\frac{4 U L}{\epsilon \min(b,h)} + 2 L \E[D].
\end{equation}
\end{claim}
\proof{Proof of Claim\ \ref{goodtime1}} 
Note that we may (deterministically) partition the time interval $[(1 - \epsilon) T_{\epsilon} - L , T_{\epsilon}]$ into $\lceil \frac{\epsilon T_{\epsilon}}{L} \rceil$ disjoint intervals each of length $L$, plus an additional disjoint time interval of length possibly less than $L$.  Suppose for contradiction that of these $\lceil \frac{\epsilon T_{\epsilon}}{L} \rceil$ disjoint time intervals of length $L$, there does not exist a single such interval $I$ s.t. 
\begin{equation}\label{inthas}
\E[|\hat{I}^{\pi^{*,\epsilon}}_{t - L_0} + q^{\pi^{*,\epsilon},E}_{t - L_0}|] \leq \frac{2 U L}{\epsilon \min(b,h)} + (L_0 + 1) \E[D]\ \textrm{for all}\ t \in I.
\end{equation}
In that case, by the triangle inequality, each of these $\lceil \frac{\epsilon T_{\epsilon}}{L} \rceil$ intervals contains at least one time period $t$ for which
$$
\E[|\hat{I}^{\pi^{*,\epsilon}}_{t - L_0} + q^{\pi^{*,\epsilon},E}_{t - L_0} - \sum_{i = t - L_0}^t D_i|] > \frac{2 U L}{\epsilon \min(b,h)}.$$
Hence by (\ref{Lip1}), non-negativity of costs, and the definition of $T_{\epsilon}$ we conclude that 
\begin{eqnarray*}
C(\pi^{*,\epsilon}) &>& \frac{\sum_{t = \lceil (1 - \epsilon) T_{\epsilon} \rceil - L}^{T_{\epsilon}} \E[C^{\pi^{*,\epsilon}}_t]}{T_{\epsilon}} - \frac{\epsilon}{2}
\\&>& \frac{\min(b,h) \times \frac{\epsilon T_{\epsilon}}{L} \times \frac{2 U L}{\epsilon \min(b,h)}}{T_{\epsilon}} - \frac{\epsilon}{2}\ \ \ =\ \ \ 2 U - \frac{\epsilon}{2}\ \ \ >\ \ \ \frac{3}{2} U,
\end{eqnarray*}
and thus $\textrm{OPT}(L) > \frac{3}{2} U - \frac{\epsilon}{2} > U$, a contradiction.  Let $t'$ denote the left end-point of the corresponding interval satisfying (\ref{inthas}), whose existence we have just proven by contradiction (in case of multiple such intervals, take the left-most such interval).  Now, further suppose for contradiction that there exists $k \in [1 , L - L_0 - 1]$ s.t. $\E[{\mathcal R}^{\pi^{*,\epsilon},t' -  L_0}_k] > \frac{4 U L}{\epsilon \min(b,h)} + 2 L \E[D].$  Then it would follow from the inventory update dynamics, non-negativity of order quantities, and the triangle inequality that 
$\E[|\hat{I}^{\pi^{*,\epsilon}}_{t' + k - L_0} + q^{\pi^{*,\epsilon},E}_{t' + k - L_0}|] > \frac{2 U L}{\epsilon \min(b,h)} + (L_0 + 1) \E[D],$ which would itself contradict the definition of $t'$.  Combining the above, and setting $T_{1,\epsilon} = t'$, completes the proof. $\Halmos$.
\endproof  
\subsubsection{Statement of approximate form of Theorem\ \ref{proxychain0}.}
We now formally state the aforementioned approximate version of Theorem\ \ref{proxychain0}.
\begin{lemma}\label{proxychain1}
For all $\epsilon \in \big(0, \min(\frac{1}{2},U) \big)$, one may construct  an $L - L_0 - 1$-dimensional random vector $\mathbf{\chi}^{*,\epsilon}$, an $L - L_0$-dimensional random vector $\mathbf{q}^{*,\epsilon}$, and a random variable $\mathcal{I}^{*,\epsilon}$, as well as $\lbrace D_i, i \geq 1 \rbrace$, on a common probability space s.t. the following are true.
\begin{enumerate}[(i)]
\item \label{enum1} $(\mathbf{\chi}^{*,\epsilon}, \mathbf{q}^{*,\epsilon},\mathcal{I}^{*,\epsilon})$ has finite mean, and w.p.1 $(\mathbf{\chi}^{*,\epsilon}, \mathbf{q}^{*,\epsilon})$ is non-negative.  Also, $(\mathbf{\chi}^{*,\epsilon},\mathcal{I}^{*,\epsilon})$ is independent of $\lbrace D_i, i \geq 1 \rbrace$ and $q^{*,\epsilon}_i$ is indepenent of $\lbrace D_j, j \geq i \rbrace$ for $i \in [1, L - L_0]$.
\item \label{enum1b} $\mathbf{\chi}^{*,\epsilon}_i \sim \mathbf{\chi}^{*,\epsilon}_1$ for $i \in [1 , L - L_0 - 1]$, and $q^{*,\epsilon}_i \sim q^{*,\epsilon}_1$ for $i \in [1 , L - L_0]$.
\item \label{enum1c} $\E[\mathbf{\chi}^{*,\epsilon}_1] + \E[q^{*,\epsilon}_1] = \E[D]$.
\item \label{lseq-30} For all $k \in [1,L - L_0]$,
$$\mathcal{I}^{*,\epsilon} + \sum_{i=1}^{k-1} (q^{*,\epsilon}_i + \chi^{*,\epsilon}_i - D_i) + q^{*,\epsilon}_k - \sum_{i = k}^{k + L_0} D_i
\sim \mathcal{I}^{*,\epsilon} + q^{*,\epsilon}_1 - \sum_{i=1}^{L_0 + 1} D_i.$$
\item \label{enum4} $$\textrm{OPT}(L) > c \E[q^{*,\epsilon}_1] + \bbe \left[G\left( \mathcal{I}^{*,\epsilon} + q^{*,\epsilon}_1 - \sum_{i=1}^{L_0 + 1} D_i\right)\right] - \epsilon.$$
\end{enumerate}
\end{lemma}
\subsubsection{Proof of Lemma\ \ref{proxychain1} by construction of a Markov chain with an appropriate stationary distribution.}
We now construct an appropriate Markov chain which repeatedly mimics $\pi^{*,\epsilon}$ for blocks of time of length $T_{1,\epsilon}$, and then (by a sequence of ordering decisions) brings the system back to a state distributed as its initial state (involving 
$-\sum_{i=1}^{\hat{G}} D'_{-i}$).  This is accomplished by allowing for an extra ``time-accounting" dimension in the state-space.  
While this ``clock" is between $1$ and $T_{1,\epsilon}$, the Markov chain dynamics parallel those of the inventory and pipeline 
in $\pi^{*,\epsilon}$ on $[1,T_{1,\epsilon}]$, and the clock increases by one in each period.  Whenever the clock reaches $T_{1,\epsilon}$, the Markov chain dynamics instead parallel those of a policy which first orders nothing until the truncated regular pipeline vector clears and the inventory position goes below 0, then places an expedited order to bring the inventory position to exactly 0, and finally orders nothing for an additional geometrically distributed number of time periods, where this geometric idling will preclude any pathological periodicity that might otherwise arise (ensuring existence of relevant stationary measures).  This brings the system back to a state in which the truncated regular pipeline vector is empty and the inventory position is distributed as $-\sum_{i=1}^{\hat{G}} D'_{-i}$, at which time the clock restarts to 1 and the cycle repeats, which thus yields a regenerative process.  We further note that in the associated Markov chain we will also keep track of the most recent expedited order, so that all relevant inventory and ordering costs can be expressed directly as a function of the state in the associated Markov chain.  This will allow us to apply the theory of regenerative processes to prove that the expected value of an appropriate function of the corresponding steady-state vector bounds the average cost incurred by $\pi^{*,\epsilon}$ on $[1,T_{\epsilon}]$, which itself well-approximates $\textrm{OPT}(L)$ (i.e. ensuring that Lemma 6.(\ref{enum4}) is satisfied).  Combining the above will allow us to prove that this steady-state vector satisfies the conditions of Lemma\ \ref{proxychain1}.
\proof{Proof of Lemma\ \ref{proxychain1}}
We construct an $(L - L_0 + 2)$-dimensional discrete-time Markov process $\lbrace \mathbf{Y}^{\epsilon}_t, t \geq 1 \rbrace = \lbrace (\mathbf{\chi}^{\epsilon,t}, q^{\epsilon}_t, \mathcal{I}^{\epsilon}_t, \tau^{\epsilon}_t) , t \geq 1 \rbrace$, where $\chi^{\epsilon,t}$ is an $(L - L_0 - 1)$-dimensional random vector, and $q^{\epsilon}_t, \mathcal{I}^{\epsilon}_t,$ and $\tau^{\epsilon}_t$ are random variables.  Let $\lbrace B_t, t \geq 1 \rbrace$ denote an i.i.d. sequence of Bernoulli r.v.s, each of which equals 1 w.p. $\frac{1}{2}$ and 0 w.p. $\frac{1}{2}$.  Then $\lbrace \mathbf{Y}^{\epsilon}_t, t \geq 1 \rbrace$ evolves as follows.  
\begin{itemize}
\item $\chi^{\epsilon,1} = \mathbf{0}\ \ \ ,\ \ \ {\mathcal I}^{\epsilon}_1 = - \sum_{i=1}^{\hat{G}} D'_{- i}\ \ \ ,\ \ \ q^{\epsilon}_1 = \hat{f}^{\pi^{*,\epsilon}}_{E,1}\big(\chi^{\epsilon,1},{\mathcal I}^{\epsilon}_1\big)\ \ \ ,\ \ \ 
\tau^{\epsilon}_1 = 1$.  
\end{itemize}
For $t \geq 1$, the dynamics are as follows.  
\begin{itemize}
\item $\chi^{\epsilon,t+1}_i = \chi^{\epsilon,t}_{i+1}$ for $i \in [1,L-L_0-2]$\ \ \ ,\ \ \ ${\mathcal I}^{\epsilon}_{t+1} = {\mathcal I}^{\epsilon}_t + \chi^{\epsilon,t}_1 + q^{\epsilon}_t - D_t$.
\item If $\tau^{\epsilon}_t \in [1,T_{1,\epsilon})$:
\begin{itemize}
\item $\tau^{\epsilon}_{t+1} = \tau^{\epsilon}_t + 1\ \ \ ,\ \ \ 
\chi^{\epsilon,t+1}_{L-L_0-1} = \hat{f}^{\pi^{*,\epsilon}}_{R,\tau^{\epsilon}_t}(\chi^{\epsilon,t},{\mathcal I}^{\epsilon}_t)\ \ \ ,\ \ \ q^{\epsilon}_{t+1} = \hat{f}^{\pi^{*,\epsilon}}_{E,\tau^{\epsilon}_{t+1}}(\chi^{\epsilon,t+1},{\mathcal I}^{\epsilon}_{t+1})$.  
\end{itemize}
\item If $\tau^{\epsilon}_t  = T_{1,\epsilon}$ and either $\chi^{\epsilon,t} \neq \mathbf{0}$ or ${\mathcal I}^{\epsilon}_1 > 0$: 
\begin{itemize}
\item $\chi^{\epsilon,t+1}_{L-L_0-1} = q^{\epsilon}_{t+1} = 0\ \ \ ,\ \ \ \tau^{\epsilon}_{t+1} = T_{1,\epsilon}.$
\end{itemize}
\item If $\tau^{\epsilon}_t  = T_{1,\epsilon}$ and $\chi^{\epsilon,t} = \mathbf{0}$ and ${\mathcal I}^{\epsilon}_{t+1} \leq 0$: 
\begin{itemize}
\item $\chi^{\epsilon,t+1}_{L-L_0-1} = 0\ \ \ ,\ \ \ q^{\epsilon}_{t+1} = - {\mathcal I}^{\epsilon}_{t+1}\ \ \ ,\ \ \ \tau^{\epsilon}_{t+1} = 0$.
\end{itemize}
\item If $\tau^{\epsilon}_t = 0$: 
\begin{itemize}
\item $\chi^{\epsilon,t+1}_{L-L_0-1} = 0\ \ \ \ ,\ \ \ q^{\epsilon}_{t+1} = 0\ \ \ ,\ \ \ \tau^{\epsilon}_{t+1} = B_t$. 
\end{itemize}
\end{itemize}
One may easily verify the following properties of $\lbrace \mathbf{Y}^{\epsilon}_t, t \geq 1 \rbrace$.  Let $z(x,y) \stackrel{\Delta}{=} \E[G(x + y - \sum_{i=1}^{L_0 + 1} D'_i)] + c y$, and $\hat{T}_{\epsilon}$ denote a r.v. distributed as the time between the chain's initial and second visit to a state s.t. $\tau^{\epsilon}_t = 1$.
\begin{itemize}
\item \label{firstc} It follows directly from the Markov chain dynamics that for all $t \geq 1$,
$\chi^{\epsilon,t+1}_i \sim \chi^{\epsilon,t}_{i+1}$ for $i \in [1,L-1]$.
\item \label{firstd0} Conditional on the event $\lbrace \tau^{\epsilon}_t = T_{1,\epsilon} \rbrace$, the expected number of time steps until $\tau^{\epsilon}_t = 0$ is at most $L + \frac{U_{\epsilon}}{\E[D]}$.
\item \label{firstd1} Conditional on the event $\lbrace \tau^{\epsilon}_t = 0 , \tau^{\epsilon}_{t-1} = T_{1,\epsilon}\rbrace$, it holds that (w.p.1) $\chi^{\epsilon,t+1} = \mathbf{0}$, ${\mathcal I}^{\epsilon}_{t+1} = - D_t$, and the number of time steps until $\tau^{\epsilon}_t = 1$ is distributed as $\hat{G}$.
\item \label{firstd} Conditional on the event $\lbrace \tau^{\epsilon}_t = 1 \rbrace$, it holds that $\mathbf{Y}^{\epsilon}_t \sim (\mathbf{0},0,-\sum_{i=1}^{\hat{G}} D'_{-i},1)$, and the joint distribution of $\lbrace \mathbf{Y}^{\epsilon}_i , i \in [t, t  + T_{\epsilon}-1] \rbrace$ is identical to that of $\lbrace ({\mathcal R}^{\pi^{*,\epsilon},i},q_i^{\pi^{*,\epsilon},E},\hat{I}^{\pi^{*,\epsilon}}_i,i), i \in [1,T_{1,\epsilon}] \rbrace$.
\item  \label{boundTeps}
W.p.1 $\hat{T}_{\epsilon} \geq T_{1,\epsilon}$, and $\E[\hat{T}_{\epsilon}] - T_{1,\epsilon} \leq L + \frac{U_{\epsilon}}{\E[D]} + 2.$
\item  \label{zboundeps}
$0 \leq 
\E[\sum_{t=1}^{\hat{T}_{\epsilon}} z(\mathcal{I}^{\epsilon}_t , q^{\epsilon}_t)]
-
\E[\sum_{t=1}^{T_{1,\epsilon}} z(\mathcal{I}^{\epsilon}_t , q^{\epsilon}_t)]
\leq 
U_{2,\epsilon}.$
\end{itemize}
Combining with the basic definitions associated with the theory of regenerative processes (here we refer the interested reader to \citet{Asmussen03} for an excellent overview), we conclude that 
$\lbrace \mathbf{Y}^{\epsilon}_t, t \geq 1 \rbrace$ is a discrete-time aperiodic regenerative process, with regeneration points coinciding with visits to states s.t. $\tau^{\epsilon}_t = 1$.  
Then we may conclude the following from standard results in the theory of regenerative processes (cf. \citet{Asmussen03,Thor92}).
\begin{enumerate}[(a)]
\item $\lbrace \mathbf{Y}^{\epsilon}_t, t \geq 1 \rbrace$ converges weakly (as $t \rightarrow \infty$) to a limiting random vector $\mathbf{Y}^{\epsilon}_{\infty} = (\chi^{\epsilon,\infty}, q^{\epsilon}_{\infty}, \mathcal{I}^{\epsilon}_{\infty}, \tau^{\epsilon}_{\infty})$.  
\item \label{constructa}Initializing the relevant Markov chain with initial conditions distributed as $\mathbf{Y}^{\epsilon}_{\infty}$ yields a stationary Markov process $\lbrace \overline{\mathbf{Y}}^{\epsilon}_t, t \geq 1 \rbrace = \lbrace 
(\overline{\mathbf{\chi}}^{\epsilon,t}, \overline{q^{\epsilon}}_t, \overline{\mathcal{I}}^{\epsilon}_t, \overline{\tau}^{\epsilon}_t , t \geq 1 \rbrace$.  Furthermore, it follows directly from the relevant Markov chain dynamics that we may construct $\lbrace \overline{\mathbf{Y}}^{\epsilon}_t, t \geq 1 \rbrace$ and $\lbrace D_i, i \geq 1 \rbrace$ on an appropriate probability space s.t. setting $\mathbf{\chi}^{*,\epsilon} = \overline{\mathbf{\chi}}^{\epsilon,1}, \mathcal{I}^{*,\epsilon} = \overline{\mathcal{I}}^{\epsilon}_1$, and 
$q^{*,\epsilon}_k = \overline{q}^{\epsilon}_k$ for $k \in [1,L-L_0]$ yields a random vector satsifying conditions (\ref{enum1}) - (\ref{lseq-30}) of Lemma\ \ref{proxychain1}.
\item \label{constructb}$\E[z(\mathcal{I}^{\epsilon}_{\infty} , q^{\epsilon}_{\infty})] = \frac{\E[\sum_{t=1}^{\hat{T}_{\epsilon}} z(\mathcal{I}^{\epsilon}_t , q^{\epsilon}_t)]}{\E[\hat{T}_{\epsilon}]}$.
\end{enumerate}
Further combining (\ref{constructb}) with our previous bounds for $\E[\hat{T}_{\epsilon}], \E[\sum_{t=1}^{\hat{T}_{\epsilon}} z(\mathcal{I}^{\epsilon}_t , q^{\epsilon}_t)]$, our definition of $T_{1,\epsilon}$, and some straightforward algebra (the details of which we omit) demonstrates that the same random vector exhibited in (\ref{constructa}) also satisfies condition (\ref{enum4}) of Lemma\ \ref{proxychain1}, completing the proof of the lemma.  $\Halmos$.
\endproof
\subsubsection{Proof of Theorem\ \ref{proxychain0}.}
We now complete the proof of Theorem\ \ref{proxychain0}, by taking an appropriate weak limit (as $\epsilon \downarrow 0$) of the random vectors which we have proven to satisfy the conditions of Lemma\ \ref{proxychain1}, and verifying certain interchanges of expectation and limit (in inequality form).
\proof{Proof of Theorem\ \ref{proxychain0}}
To complete the proof of Theorem\ \ref{proxychain0}, we now prove that the sequence of random vectors $\lbrace (\mathbf{\chi}^{*,\frac{1}{n}}, \mathbf{q}^{*,\frac{1}{n}}, \mathcal{I}^{*,\frac{1}{n}}) , n \geq 2 + \frac{1}{U} \rbrace$ is tight.  It follows from Lemma\ \ref{proxychain1}.(\ref{enum1c}) and (\ref{enum4}), non-negativity, the triangle inequality, the fact that $\textrm{OPT}(L) \leq U$, and (\ref{Lip1}) that for all $n \geq 2 + \frac{1}{U}$, 
\begin{equation}\label{boundIn}
\E[|\mathcal{I}^{*,\frac{1}{n}}|] \leq \frac{2 U}{\min(b,h)} + (L_0 + 1)\E[D].
\end{equation}
Combining with Lemma\ \ref{proxychain1}.(\ref{enum1b}) - (\ref{enum1c}) and non-negativity, we conclude the desired tightness, and hence existence of at least one subsequential limit (cf. \cite{Billingsley99}) $(\mathbf{\chi}^{*,\infty}, \mathbf{q}^{*,\infty}, \mathcal{I}^{*,\infty})$.  Let $\lbrace n_i, i \geq 1 \rbrace$ denote any fixed subsequence along which the sequence of measures converges to this limit, s.t. $n_1 \geq 2 + \frac{1}{U}$.  That this weak limit satisfies Theorem\ \ref{proxychain0}.(\ref{enum0}) - (\ref{lseq-300}) follows from the definition of weak convergence.  However, it will require somewhat subtle reasoning to prove (\ref{enum0mean1}) - (\ref{enum0mean3}), since e.g. dominated convergence does not necessarily hold and thus one must take care when interchanging limit and expectation.  Note that by the Skorohod representation theorem and continuous mapping theorem, we may construct $\lbrace (\mathbf{\chi}^{*,\frac{1}{n_i}}, \mathbf{q}^{*,\frac{1}{n_i}}, |\mathcal{I}^{*,\frac{1}{n_i}}|) , i \geq 1 \rbrace$ and $(\mathbf{\chi}^{*,\infty}, \mathbf{q}^{*,\infty}, |\mathcal{I}^{*,\infty}|)$ on a common probability space s.t. the corresponding sequence of random vectors converges almost surely (i.e. not only in distribution) to $(\mathbf{\chi}^{*,\infty}, \mathbf{q}^{*,\infty}, |\mathcal{I}^{*,\infty}|)$.  As all associated r.v.s are non-negative, we may apply Fatou's lemma to conclude that $\E[\mathbf{\chi}^{*,\infty}_1] \leq \liminf_{i \rightarrow \infty} \E[\chi^{*,\frac{1}{n_i}}], \E[q^{*,\infty}_1] \leq \liminf_{i \rightarrow \infty} \E[q^{*,\frac{1}{n_i}}_1]$, and $\E[|\mathcal{I}^{*,\infty}|] \leq \liminf_{i \rightarrow \infty} \E[|\mathcal{I}^{*,\frac{1}{n_i}}|]$.  Combining with Lemma\ \ref{proxychain1}.(\ref{enum1}) and (\ref{enum1c}), as well as (\ref{boundIn}), then completes the proof of Theorem\ \ref{proxychain0}.(\ref{enum0mean1}).  Combining with the already proven Theorem\ \ref{proxychain0}.(\ref{enum0b}) and (\ref{lseq-300}), with $k = 2$, yields Theorem\ \ref{proxychain0}.(\ref{enum0mean2}).  Finally, we prove that the corresponding vector also satisfies Theorem\ \ref{proxychain0}.(\ref{enum0mean3}).  Let $Z_n \stackrel{\Delta}{=} c q^{*,\frac{1}{n}}_1 + G\left( \mathcal{I}^{*,\frac{1}{n}} + q^{*,\frac{1}{n}}_1 - \sum_{i=1}^{L_0 + 1} D_i\right)$, and $Z_{\infty} \stackrel{\Delta}{=}  c q^{*,\infty}_1 + G\left( \mathcal{I}^{*,\infty} + q^{*,\infty}_1 - \sum_{i=1}^{L_0 + 1} D_i\right)$.  The already proven weak convergence, and continuous mapping theorem, implies that $\lbrace Z_{n_i}, i \geq 1 \rbrace$ converges weakly to $Z_{\infty}$.  It follows from the Skorohod representation theorem (cf. \cite{Billingsley99}) that we may construct $\lbrace Z_{n_i}, i \geq 1 \rbrace$ and $Z_{\infty}$ on a common probability space so that this convergence holds almost surely (as opposed to only in distribution).  Applying non-negativity and Fatou's lemma, we conclude that on this probability space, $\E[\liminf_{i \rightarrow \infty} Z_{n_i}] \leq \liminf_{i \rightarrow \infty} E[Z_{n_i}]$, and hence (combining with the stated almost sure convergence) $\E[Z_{\infty}] \leq \liminf_{i \rightarrow \infty} E[Z_{n_i}]$.  Combining with Lemma\ \ref{proxychain1}.(\ref{enum4}), which implies that $\textrm{OPT}(L) > \E[Z_{n_i}] - \frac{1}{n_i}$ for all $i \geq 1$, and the already proven Theorem\ \ref{proxychain0}.(\ref{enum0mean2}), completes the proof.  $\Halmos.$
\endproof

\subsection{Proof of Lemma\ \ref{limequal2}}
In preparation for bounding $V_{\alpha}^{\infty}(r,x) - V_{\alpha}^n(r,x)$, we first bound the optimal value, and set of minimizers, of $V_{\alpha}^n(r,x)$, uniformly in $n$.  
For $\alpha \in (0,1)$ and $r \in \bbr$, let $\overline{S}_{\alpha}^n(r)$ denote the supremum of the set of minimizers (with respect to $x$) of $V_{\alpha}^{n}(r,x)$, where we note that a straightforward contradiction demonstrates that $\overline{S}_{\alpha}^n(r) \in (-\infty, \infty)$ for each $\alpha,n,r$; and it follows from Lemma\ \ref{use1} that $V_{\alpha}^{n}(r,-\infty) = V_{\alpha}^{n}\big(r,\overline{S}_{\alpha}^n(r)\big)$.  Then we prove the following uniform bounds.  

\begin{lemma}\label{use2}
\begin{enumerate}
\item \label{use21} For $\alpha \in (0,1)$ and $r, x \in \bbr$, it holds that 
$$\sup_{n \geq 1} V_{\alpha}^{n}(r,x) < 2 (L_0 + 1) \max(b,h) (|x| + |r| + \E[D]) (1 - \alpha)^{-2}.$$  
\item \label{use22} For $\alpha \in (0,1)$ and $n \geq 1$, it holds that $|\overline{S}_{\alpha}^n(r)| < \overline{S}_{\alpha}(r)$. 
\item \label{use23} For all $y \notin[- \overline{S}_{\alpha}(r), \overline{S}_{\alpha}(r)]$ and $n \geq 1$,
$$
\bbe\left[G\big(y - \sum_{k=n}^{L_0 + n} (D_k - r) \big)\right] + \alpha \bbe\left[V_{\alpha}^{n-1}\big(r, y - (D_{L_0 + n} - r)\big)\right] > V_{\alpha}^n\big(r,\overline{S}_{\alpha}^n(r)\big) +  (L_0 + 1) \max(b,h) \E[D].$$
\item \label{use24} For all $L > L_0 + 1$, 
$$\textrm{OPT}(L) \geq c(\bbe[D] - r_L) + (1 - \alpha) V_{\alpha}^{L - L_0}\big(r_L, - \overline{S}_{\alpha}(r_L) \big).$$
\end{enumerate}
\end{lemma}
\proof{Proof of Lemma\ \ref{use2}}
By evaluating the policy which never orders, we conclude that for all $\alpha \in (0,1)$, $r,x \in \bbr$, $\sup_{n \geq 1} V_{\alpha}^{n}(r,x)$ is at most $\bbe\left[\sum_{i=1}^{\infty}\alpha^{i-1}G\big(x - \sum_{j=1}^i (D_j - r) - \sum_{k = i+1}^{L_0 + i} (D_k - r)\big)\right]$, which by (\ref{Lip1}) is itself bounded by
$$\max(b,h) (|x| + |r| + \E[D]) \sum_{i=1}^{\infty} (i + L_0) \alpha^{i-1}\ \ \ <\ \ \ 2 (L_0 + 1) \max(b,h) (|x| + |r| + \E[D]) (1 - \alpha)^{-2}.
$$
The remainder of the lemma follows from (\ref{Lip1}), Lemmas\ \ref{reduce0} and \ref{use1}, and a straightforward calculation and argument by contradiction, and we omit the details.  Combining the above completes the proof.  $\Halmos$
\endproof

With Lemma\ \ref{use2} in hand, we now complete the proof of Lemma\ \ref{limequal2}.

\proof{Proof of Lemma\ \ref{limequal2}}
We first demonstrate that $V_{\alpha}^{\infty}(r,x) = \lim_{n \rightarrow \infty} V_{\alpha}^n(r,x)$, and complete the proof of (\ref{limequal2a}).  The existence of the corresponding limit follows from the monotonicity (in $n$) guaranteed by Lemma\ \ref{use1}.
That $V_{\alpha}^{\infty}(r,x) \geq \lim_{n \rightarrow \infty} V_{\alpha}^n(r,x)$ for all $\alpha \in (0,1)$ and $r,x \in \bbr$ follows immediately from the definitions of the associated optimization problems.  To prove the other direction, as well as (\ref{limequal2a}), we note that for any fixed $n \geq 1$, it follows from the convexity ensured by Lemma\ \ref{use1} that there exists an optimal policy $\pi$ for the problem stated in the r.h.s. of (\ref{ds-finite}) of base-stock form, with order-up-to levels $C_1,\ldots,C_n$ (i.e. order up to level $C_i$ in period $i$ if the pre-order inventory level is below $C_i$, otherwise order nothing).  Furthermore, it follows from Lemma\ \ref{use2} that $\max_{i = 1,\ldots,n} |C_i| \leq \overline{S}_{\alpha}(r)$.  Now, consider the policy $\pi'$ $\big($ for the problem stated in the r.h.s. of (\ref{ds-infinite}) $\big)$ that orders up to level $C_i$ in period $i$ if the pre-order inventory position is below $C_i$ and otherwise orders nothing in periods $i=1,\ldots,n$; and orders nothing in all remaining periods, irregardless of the inventory level.  It follows from a straightforward bounding argument that under policy $\pi'$, w.p.1 the absolute value of the post-ordering inventory position in period $i$ is at most $|x| + \overline{S}_{\alpha}(r) + (i-1) |r| + \sum_{k=1}^{i-1} D_k$.  Thus by the dynamics of the underlying inventory problem and (\ref{Lip1}), it follows that for all $i \geq n + 1$, $C^{\pi'}_i \leq \max(b,h)\big( |x| + \overline{S}_{\alpha}(r) + (i + L_0) |r| + (i + L_0) \E[D] \big)$.  Thus since (by construction) $C^{\pi'}_i = C^{\pi}_i$ for $i \in [1,n]$, it follows from definitions and straightforward algebra that
\begin{eqnarray*}
\bbe\left[\sum_{ i = 1}^{\infty} \alpha^{i-1} C^{\pi'}_i(r,x)\right]- V_{\alpha}^n(r,x) &\leq& 
\max(b,h) (\overline{S}_{\alpha}(r) + |x| + |r| + \E[D]) \sum_{i=n + 1}^{\infty} (i + L_0) \alpha^{i-1} 
\\&\leq& \max(b,h) (\overline{S}_{\alpha}(r) + |x| + |r| + \E[D]) (1 + L_0 + n) (1 - \alpha)^{-2} \alpha^n.
\end{eqnarray*}
This completes the proof of (\ref{limequal2a}), and letting $n \rightarrow \infty$ completes the proof that $V_{\alpha}^{\infty}(r,x) = \lim_{n \rightarrow \infty} V_{\alpha}^n(r,x)$.  Combining with Lemmas\ \ref{use1} and\ \ref{use2}, the fact that convexity and monotonicity are preserved under limits, and a straightforward contradiction argument completes the proof of all parts of the lemma regarding properties of $V_{\alpha}^{\infty}(r,x)$ and the associated optimization problems and optimal policies.  Finally, we complete the proof of (\ref{limequal2b}).  It follows from Lemmas\ \ref{use1} and\ \ref{use2}, the already proven parts of Lemma\ \ref{limequal2}, and the fact that Theorem\ \ref{proxychain0} ensures $r_L \in \big[0, \E[D] \big]$ that $\textrm{OPT}(L) - c(\bbe[D] - r_L)$ is at least
$$(1 - \alpha) \bigg( V_{\alpha}^{\infty}\big(r_L, S_{\alpha}^{\infty}(r_L) \big) - 2 \max(b,h) \big(\overline{S}_{\alpha}(\E[D]) + \E[D]\big) (1 + L) (1 - \alpha)^{-2} \alpha^{L - L_0} \bigg).$$
Combining with some straightforward algebra, the definition of $U_0$, and the already proven parts of Lemma\ \ref{limequal2} completes the proof.  $\Halmos$.
\endproof

\subsection{Proof of Corollary\ \ref{lbound1}}
Before proving Corollary\ \ref{lbound1}, we will need a preliminary result which demonstrates that if $r$ is ``very close" to $\E[D]$, then $M^r_i$ is ``very large" for an appropriate range of $i$.  We will then use this result to show that $r_L$ cannot be ``too close" to $\E[D]$ by deriving a contradiction, showing that in this case the optimal value would be strictly greater than $U$, which is impossible.  

\begin{lemma}\label{zreviewa}
If there exists $\epsilon \in \big[0, \E[D] - Q_0\big]$ s.t. $r \in \big(\E[D] - \epsilon, \E[D]\big]$, then for all  $i,j \in \big[400 \hat{p}_0^{-2}, (\hat{p}_0 \eta_0 \epsilon^{-1})^2 \big]$ s.t. $i \geq j$, 
$$M^r_i - M^r_j \geq \frac{1}{5} \hat{p}_0 \eta_0  \big( i^{\frac{1}{2}} - j^{\frac{1}{2}} \big) - (i - j) \epsilon - 2 \eta_0 \big(\log(\frac{i}{j}) + 2\big).
$$
\end{lemma}
\proof{Proof of Lemma\ \ref{zreviewa}}
We note that the result would follow from well-known weak-convergence results under additional assumptions on $D$ (e.g. finite variance, cf. \citet{Erdos46}), but to avoid unnecessary assumptions (and for completeness) we provide a proof from first principles.  
Let us fix any $\epsilon \in \big[0, \E[D] - Q_0\big]$, $r \in \big(\E[D] - \epsilon, \E[D]\big]$, and $k \in \big[400 \hat{p}_0^{-2}, (\hat{p}_0 \eta_0 \epsilon^{-1})^2 \big]$ (supposing this interval is non-empty, i.e. $\epsilon \leq \frac{\hat{p}_0^2 \eta_0}{20}$).  Let $\lbrace A^{+,r}_i, i \geq 1 \rbrace$ denote an i.i.d. sequence of r.v.s distributed as $r - D$ conditioned on the event $\lbrace r > D \rbrace$, and $\lbrace A^{-,r}_i, i \geq 1 \rbrace$ denote an i.i.d. sequence of r.v.s distributed as $D - r$ conditioned on the event $\lbrace D \geq r \rbrace$.  Let $B^r_k$ denote a binomially distributed r.v. with parameters $k, \rho_r \stackrel{\Delta}{=} \p(\lbrace r > D \rbrace)$, independent of $\lbrace A^{+,r}_i, i \geq 1 \rbrace$ and $\lbrace A^{-,r}_i, i \geq 1 \rbrace$.  It follows from definitions and the constraints on $\epsilon$ and $r$ that
$$\rho_r \in [\frac{1}{2} p_0, p_0].
$$
Note that for $k \geq 1$, we may construct $W^r_k$ on an appropriate probability space s.t. $W^r_k = \sum_{i=1}^{B^r_k} A^{+,r}_i - \sum_{i=1}^{k - B^r_k} A^{-,r}_i$, in which case (by non-negativity) $\E[\max(0,W^r_k)]$ is at least
$$\E\bigg[\sum_{i=1}^{B^r_k} A^{+,r}_i - \sum_{i=1}^{k - B^r_k} A^{-,r}_i \bigg| \bigg\lbrace B^r_k \geq \rho_r k + \big( \rho_r (1 - \rho_r) k \big)^{\frac{1}{2}} \bigg\rbrace \bigg]
\p\bigg( \frac{B^r_k - \rho_r k}{\big(\rho_r (1 - \rho_r) k\big)^{\frac{1}{2}}} \geq 1 \bigg).$$
Furthermore, since $\rho_r \E[A^{+,r}_1] = (1 - \rho_r)\E[A^{-,r}_1] - (\E[D] - r)$, it follows from non-negativity and independence that
\begin{eqnarray*}
\ &\ &\ \E\bigg[\sum_{i=1}^{B^r_k} A^{+,r}_i - \sum_{i=1}^{k - B^r_k} A^{-,r}_i \bigg| \bigg\lbrace B^r_k \geq \rho_r k + \big( \rho_r (1 - \rho_r) k \big)^{\frac{1}{2}} \bigg\rbrace \bigg]
\\&\ &\ \ \ \ \ \ \geq\ \ \ \big( \rho_r k + \big( \rho_r (1 - \rho_r) k \big)^{\frac{1}{2}}\big) \E[A^{+,r}_1] - \big( (1 - \rho_r) k - \big( \rho_r (1 - \rho_r) k \big)^{\frac{1}{2}}\big) \E[A^{-,r}_1]
\\&\ &\ \ \ \ \ \ =\ \ \ \big( \rho_r (1 - \rho_r) k \big)^{\frac{1}{2}} ( \E[A^{+,r}_1] + \E[A^{-,r}_1] ) - k (\E[D] - r).
\end{eqnarray*}
Let $N(0,1)$ denote a standard normal r.v.  By the celebrated Berry-Esseen Theorem (cf. \cite{Korolev10}),
$$
\bigg|\p\bigg( \frac{B^r_k - \rho_r k}{\big(\rho_r (1 - \rho_r) k\big)^{\frac{1}{2}}} \geq 1 \bigg) - \p\big( N(0,1) \geq 1 \big) \bigg| \leq 2 \big(\rho_r(1 - \rho_r) k \big)^{-\frac{1}{2}}.
$$
It is easily verified from definitions that 
$$
\E[A^{+,r}_1] + \E[A^{-,r}_1]  \geq \eta_0\ \ \ ,\ \ \ \big( \rho_r (1 - \rho_r) \big)^{\frac{1}{2}} \geq \hat{p}_0\ \ \ ,\ \ \ \p\big( N(0,1) \geq 1 \big) \geq \frac{1}{10}.
$$
Thus combining the above, we conclude that 
\begin{equation}\label{part1a}
\E\bigg[\sum_{i=1}^{B^r_k} A^{+,r}_i - \sum_{i=1}^{k - B^r_k} A^{-,r}_i \bigg| \bigg\lbrace B^r_k \geq \rho_r k + \big( \rho_r (1 - \rho_r) k \big)^{\frac{1}{2}} \bigg\rbrace \bigg] \geq \hat{p}_0 \eta_0 k^{\frac{1}{2}} - \epsilon k,
\end{equation}
and
\begin{equation}\label{part1b}
\p\bigg( \frac{B^r_k - \rho_r k}{\big(\rho_r (1 - \rho_r) k\big)^{\frac{1}{2}}} \geq 1 \bigg) \geq \frac{1}{10} - 2 \hat{p}_0^{-1} k^{-\frac{1}{2}}.
\end{equation}
As our assumptions on $\epsilon,r,k$ ensure that the r.h.s. of both (\ref{part1a}) and (\ref{part1b}) are non-negative, we conclude that 
$\E[\max(0,W^r_k)]$ is at least $\big(\hat{p}_0 \eta_0 k^{\frac{1}{2}} - \epsilon k\big)\big(\frac{1}{10} - 2 \hat{p}_0^{-1} k^{-\frac{1}{2}}\big)$, which is itself at least $\frac{1}{10} \hat{p}_0 \eta_0 k^{\frac{1}{2}} - k \epsilon - 2 \eta_0$.  Thus by Lemma\ \ref{zreview},
\begin{eqnarray*}
M_i - M_j &=& \sum_{l=j}^{i-1} l^{-1} \E\big[\max(0,W^r_l)\big]
\\&\geq& \frac{1}{10} \hat{p}_0 \eta_0 \sum_{l = j}^{i-1} l^{- \frac{1}{2}} - (i - j) \epsilon - 2 \eta_0 \sum_{l=j}^{i-1} l^{-1}
\\&\geq& \frac{1}{10} \hat{p}_0 \eta_0 \int_j^{i} x^{- \frac{1}{2}} dx - (i - j) \epsilon - 2 \eta_0 \big(\log(\frac{i}{j}) + 2\big)
\\&=& \frac{1}{5} \hat{p}_0 \eta_0  \big( i^{\frac{1}{2}} - j^{\frac{1}{2}} \big) - (i - j) \epsilon - 2 \eta_0 \big(\log(\frac{i}{j}) + 2\big),
\end{eqnarray*}
where we have used the well-known fact that for all $n \geq 1$, $\log(n) \leq \sum_{l=1}^n l^{-1} \leq \log(n) + 2$.  Combining the above completes the proof.  $\Halmos$
\endproof
Before completing the proof of Corollary\ \ref{lbound1}, it will be useful to collect a few additional auxiliary bounds.  For $\alpha \in (0,1)$, let $G_{\alpha}$ denote a geometrically distributed r.v. with success probability $1 - \alpha$, i.e. $\p(G_{\alpha} = k) = (1 - \alpha) \alpha^{k-1}$ for $k \geq 1$, independent of $\lbrace Z^{r_L}_k, k \geq 1 \rbrace$, and $m_{\alpha} \stackrel{\Delta}{=} \lceil - \frac{1}{\log_2(\alpha)} \rceil$ denote a median of $G_{\alpha}$.  Note that the memoryless property implies $\p\big(G_{\alpha} \geq 2 m_{\alpha}\big) \geq \frac{1}{4}$.  Let $\underline{\xi}_0 \stackrel{\Delta}{=} 2^{-\frac{\hat{p}^2_0}{400}}$, and $\overline{\xi}_0 \stackrel{\Delta}{=} 2^{- \frac{4}{\hat{p}^2_0 \eta_0^2} \epsilon_0^2}$.  
\begin{lemma}\label{almostcor1}
\ \begin{enumerate}[(i)]
\item $.998 < \underline{\xi}_0 \leq 1 - \epsilon_0 \leq \overline{\xi}_0 < 1$. \label{prover1}
\item $L \geq \epsilon_0^{-2}$ implies $\epsilon_0^{-3} L \exp(- \epsilon_0 L) \leq 25$. \label{prover2}
\item $\alpha \in [\underline{\xi}_0, \overline{\xi}_0]$ implies: \label{malpha1}
	\begin{itemize} 
		\item $m_{\alpha}\ ,\ 2 m_{\alpha}\ \in\ \big[400 \hat{p}_0^{-2}, (\hat{p}_0 \eta_0 \epsilon_0^{-1})^2 \big]$;
		 \item $\frac{1}{4} (1 - \alpha)^{-1} \leq m_{\alpha} \leq 4 (1 - \alpha)^{-1}.$ 
	\end{itemize} 
\end{enumerate}
\end{lemma}
\proof{Proof of Lemma\ \ref{almostcor1}}
\ \\\ref{prover1}: That $\underline{\xi}_0 > .998$ follows from the fact that $\hat{p}_0 < 1$, and thus $2^{-\frac{\hat{p}^2_0}{400}} > 2^{-\frac{1}{400}} \geq .998$.  That $\underline{\xi}_0 \leq 1 - \epsilon_0$ follows from the fact that (by definition) $\epsilon_0 \leq 1 - \underline{\xi}_0$.  We now prove that $1 - \epsilon_0 \leq \overline{\xi}_0$.  By definition $\epsilon_0 \leq \frac{1}{4}(\eta_0 \hat{p}_0)^2$, which implies that $\frac{4}{\hat{p}^2_0 \eta_0^2} \epsilon_0^2 \leq \epsilon_0$.  Combining with the exponential inequality and the fact that $\log(2) < 1$, we conclude that 
$$\overline{\xi}_0 \geq 2^{- \epsilon_0} \geq 1 - \log(2) \epsilon_0 \geq 1 - \epsilon_0,$$
completing the proof.  As trivially $\overline{\xi}_0 < 1$, this completes the demonstration.
\\\\\ref{prover2}: It is easily verified that $\zeta_1(L) \stackrel{\Delta}{=} L\exp(- \epsilon_0 L)$ is decreasing in $L$ on $[\epsilon_0^{-1},\infty)$.  Thus $L \geq \epsilon_0^{-2}$ implies 
$$
\epsilon_0^{-3} L \exp(- \epsilon_0 L) \leq \epsilon_0^{-3} \zeta_1(\epsilon_0^{-2}) = \epsilon_0^{-5} \exp( -\epsilon_0^{-1} ).
$$
As $\zeta_2(\epsilon_0) \stackrel{\Delta}{=} \epsilon_0^{-5} \exp(- \epsilon_0^{-1})$ is increasing in $\epsilon_0$ on $(0,\frac{1}{5})$, and by definition $\epsilon_0 < \frac{1}{5}$, it follows that $\zeta_2(\epsilon_0) \leq \zeta_2(\frac{1}{5}) < 25$.  Combining the above completes the proof.
\\\\\ref{malpha1}: The first assertion follows immediately from the definitions of $\underline{\xi}_0, \overline{\xi}_0$, and $m_{\alpha}$,  and a straightforward calculation.  To prove the second assertion, note that due to (\ref{prover1}), 
$\alpha \in [\underline{\xi}_0, \overline{\xi}_0]$ implies $\alpha \in (.998 , 1)$.  It follows from a straightforward Taylor expansion of the logarithm function that $\alpha \in (.998 , 1)$ implies $ - 2(1 - \alpha) \leq \log_2(\alpha) \leq -(1 - \alpha)$, and thus 
$\frac{1}{2} (1 - \alpha)^{-1} \leq - \frac{1}{\log_2(\alpha)} \leq (1 - \alpha)^{-1}$.  Noting that $\alpha \in (.998 , 1)$ implies 
$\lceil  (1 - \alpha)^{-1} \rceil \leq 4 (1 - \alpha)^{-1}$ completes the proof.
\endproof
With Lemmas\ \ref{zreviewa} and\ \ref{almostcor1} in hand, we now complete the proof of Corollary\ \ref{lbound1}.
\proof{Proof of Corollary\ \ref{lbound1}}
Suppose for contradiction that for some $L > \epsilon_0^{-2} + L_0 + 1$, it holds that $r_L > \E[D] - \epsilon_0$.  In this case, it follows from Corollary\ \ref{wrapit0}, (\ref{Lip1}), and Jensen's inequality that for all $\alpha \in (\frac{1}{2},1)$,
\begin{equation}\label{geonews}
\textrm{OPT}(L) \geq \min(b,h) \inf_{S \in \bbr} \sum_{k=1}^{\infty} (1 - \alpha) \alpha^{k-1} \big|S + M^{r_L}_k \big| - U_0 2^{L_0} (1 - \alpha)^{-3} L \alpha^L.
\end{equation}
Note that we may interpret the r.h.s. of (\ref{geonews}) as an appropriate single-stage newsvendor problem (with ordering level $S$ and demand distributed as $M^{r_L}_{G_{\alpha}}$).  
We conclude from Lemmas\ \ref{zreviewa} and\ \ref{almostcor1}, well-known results for the newsvendor problem (cf. \citet{Zipkin00}), and the memoryless property that for all $\alpha \in [\underline{\xi}_0, \overline{\xi}_0]$,
\begin{eqnarray*}
\textrm{OPT}(L) &\geq& \min(b,h) \E\big[\big|M^{r_L}_{m_{\alpha}} - M^{r_L}_{G_{\alpha}}\big|\big] - U_0 2^{L_0} (1 - \alpha)^{-3} L \alpha^L
\\&\geq& \frac{1}{4} \min(b,h) \big(M^{r_L}_{2 m_{\alpha}} - M^{r_L}_{m_{\alpha}} \big) - U_0 2^{L_0} (1 - \alpha)^{-3} L \alpha^L
\\&\geq& \frac{1}{4} \min(b,h) \bigg( \frac{1}{5} \hat{p}_0 \eta_0  \big( (2 m_{\alpha})^{\frac{1}{2}} - m_{\alpha}^{\frac{1}{2}} \big) - \epsilon_0 m_{\alpha} - 2 \eta_0 \big(\log(2) + 2\big) \bigg) - U_0 2^{L_0} (1 - \alpha)^{-3} L \alpha^L
\\&\geq& \frac{1}{20} \min(b,h) \hat{p}_0 \eta_0 (2^{\frac{1}{2}} - 1) m^{\frac{1}{2}}_{\alpha} - \epsilon_0 m_{\alpha} - 6 \eta_0 - U_0 2^{L_0} (1 - \alpha)^{-3} L \alpha^L
\\&\geq& c_0 (1 - \alpha)^{-\frac{1}{2}} - 4 \epsilon_0 (1 - \alpha)^{-1} - 6 \eta_0 - U_0 2^{L_0} (1 - \alpha)^{-3} L \alpha^L.
\end{eqnarray*}
Setting $\alpha = 1 - \epsilon_0$, and combining the above with Lemma\ \ref{almostcor1}.(\ref{prover1}) and the fact that $1 - \epsilon_0 \leq \exp(- \epsilon_0)$, we conclude that 
$$\textrm{OPT}(L) \geq c_0 \epsilon_0^{-\frac{1}{2}} - U_0 2^{L_0} \epsilon_0^{-3} L \exp(- \epsilon_0 L) - 6 (\eta_0 + 1).
$$
Applying Lemma\ \ref{almostcor1}.(\ref{prover2}) and the fact that $\textrm{OPT}(L) \leq U$, along with a straightforward contradiction argument, completes the proof.  $\Halmos$.
\endproof

\end{document}